\documentclass[final]{colt2024}
\usepackage{times}
\usepackage{mathabx}

\renewcommand{\leq}{\leqslant}

\renewcommand{\geq}{\geqslant}
\renewcommand{\succeq}{\succcurlyeq}

\def\argmin{\operatornamewithlimits{argmin}}

\newcommand{\eps}{\varepsilon}

\newcommand{\Be}{\mathrm{Be}}
\newcommand{\dd}{\mathrm{d}}
\newcommand{\Fr}{\mathrm{F}}
\newcommand{\rank}{\mathrm{rank}}
\newcommand{\Tr}{\mathrm{Tr}}
\newcommand{\Var}{\mathrm{Var}}
\renewcommand{\vec}{\mathrm{\mathbf{vec}}}

\newcommand{\ttr}{\mathtt{r}}
\newcommand{\ttZ}{\mathtt Z}

\newcommand{\E}{\mathbb E}
\newcommand{\p}{\mathbb P}
\newcommand{\R}{\mathbb R}
\newcommand{\X}{\mathbb X}
\newcommand{\Y}{\mathbb Y}

\newcommand{\I}{\mathcal I}
\newcommand{\KL}{{\mathcal{KL}}}
\newcommand{\N}{\mathcal N}
\newcommand{\cO}{\mathcal O}

\newcommand{\ro}{\mathcal R} 
\newcommand{\cX}{{\mathcal{X}}}

\newcommand{\bu}{{\mathbf u}}
\newcommand{\bv}{{\mathbf v}}
\newcommand{\bX}{{\mathbf X}}
\newcommand{\bgamma}{{\boldsymbol \gamma}}
\newcommand{\btheta}{{\boldsymbol \theta}}

\newcommand{\bzeta}{{\boldsymbol \zeta}}
\newcommand{\bzero}{{\boldsymbol 0}}

\newtheorem{Th}{Theorem}[section]
\newtheorem{Lem}[Th]{Lemma}
\newtheorem{Def}[Th]{Definition}
\newtheorem{Prop}[Th]{Proposition}

\newtheorem{Rem}[Th]{Remark}
\newtheorem{Ex}[Th]{Example}
\newtheorem{As}[Th]{Assumption}

\title[Dimension-free Structured Covariance Estimation]{Dimension-free Structured Covariance Estimation}

\coltauthor{%
 \Name{Nikita Puchkin} \Email{npuchkin@hse.ru}\\
 \addr HSE University and IITP RAS, Russian Federation
 \AND
 \Name{Maxim Rakhuba} \Email{mrakhuba@hse.ru}\\
 \addr HSE University, Russian Federation%
}

\begin{document}

\maketitle

\begin{abstract}
    Given a sample of i.i.d. high-dimensional centered random vectors, we consider a problem of estimation of their covariance matrix $\Sigma$ with an additional assumption that $\Sigma$ can be represented as a sum of a few Kronecker products of smaller matrices. Under mild conditions, we derive the first non-asymptotic dimension-free high-probability bound on the Frobenius distance between $\Sigma$ and a widely used penalized permuted least squares estimate. Because of the hidden structure, the established rate of convergence is faster than in the standard covariance estimation problem.
\end{abstract}


\section{Introduction}
\label{sec:introduction}

Let $\bX, \bX_1, \dots, \bX_n \in \R^d$ be i.i.d. centered random vectors. We are interested in estimation of the covariance matrix $\Sigma = \E \bX \bX^\top$ from observations $\bX_1, \dots, \bX_n$. This classical problem is of a great significance and it arises in various domains such as wireless communications \citep{werner07}, economics and finance \citep{bai11}, biology and medicine (for example, functional MRI \citep{derado10}, genomics \citep{xie03}, MEG/EEG \citep{de04}). A standard approach to this problem is to consider the sample covariance
\[
    \widehat \Sigma = \frac1{n} \sum\limits_{i = 1}^n \bX_i \bX_i^\top.
\]
Unfortunately, if the number of available samples $n$ is much less than the ambient dimension $d$ (which is often the case in practical scenarios), the estimate $\widehat \Sigma$ suffers from the curse of dimensionality. To overcome this issue, researchers have explored leveraging structural properties of covariance matrices to reduce their total number of degrees of freedom. This includes assumptions on the covariance matrix like sparsity \citep{banerjee08,ravikumar11}, low-rank \citep{fan08}, Toeplitz \citep{fuhrmann91} and bandable model \citep{bickel08a, bickel08b}. A reader is referred to the survey \citep{cai16} for further examples. In situations when data can be represented as tensors, a popular modelling assumption is based on expansion of $\Sigma$ into a Kronecker product of two smaller matrices \citep{wang22}:
\begin{equation}
    \label{eq:sigma_kronecker_simple}
    \Sigma = \Phi \otimes \Psi,
    \quad
    \text{where $\Phi \in \R^{p \times p}$, $\Psi \in \R^{q \times q}$, and $pq = d$.}
\end{equation}
Each of the matrices $\Phi, \Psi$ corresponds to a certain index (mode) of the multivariate array and can be interpreted as a mode-wise covariance. For instance, when researchers deal with spatio-temporal data, they use Kronecker product to separate spatial and temporal dimensions. Despite its simplicity, the Kronecker product model \eqref{eq:sigma_kronecker_simple} appears to be useful in various applications \citep{werner07, werner08, allen10, guggenberger23}.
In the present paper, we are concerned with a more expressive model that approximates the covariance matrix as a sum of Kronecker products:
\begin{equation}
    \label{eq:sigma_kronecker}
    \Sigma = \Phi_1 \otimes \Psi_1 + \dots + \Phi_K \otimes \Psi_K,
\end{equation}
where $\Phi_1, \Psi_1, \dots, \Phi_K, \Psi_K$ are symmetric positive semidefinite matrices, such that $\Phi_j \in \R^{p \times p}$, $\Psi_j \in \R^{q \times q}$ for all $j \in \{1, \dots, K\}$ and $pq = d$. In this case, one says that $\Sigma$ has a \emph{Kronecker rank} equal to $K$.
By summing multiple Kronecker products, the model can capture complex interactions and dependencies that may not be adequately represented by a single Kronecker product. We would like to note that any symmetric matrix of size $(pq \times pq)$ admits an expansion of the form \eqref{eq:sigma_kronecker} with $K = \min\{p^2, q^2\}$. However, even with a small number of Kronecker product terms, the model \eqref{eq:sigma_kronecker} has shown promising results in applications such as clutter covariance estimation for airborne STAP \citep{sun18}, MEG/EEG \citep{de02,de04,bijma05}, video modelling \citep{greenewald13}, and SAR data analysis \citep{rucci10}. Notably, accurate approximations with a small number of Kronecker product terms have been demonstrated for the Mat\'ern family of covariances \citep{litvinenko19}. 

The rates of convergence in the covariance estimation problem with $\Sigma$ of the form \eqref{eq:sigma_kronecker_simple} or \eqref{eq:sigma_kronecker} were studied in numerous papers (see, for example, \cite{werner08, tsiligkaridis13, zhou14, leng18, masak22}). However, the existing state-of-the-art results are either dimension-dependent or based on unsuitable for high-dimensional setting assumptions. At the same time, the recent advances in covariance estimation \citep{bunea15, puchkin23, zhivotovskiy24} allow to establish non-asymptotic high-probability bounds on the Frobenius and operator norm of $(\widehat \Sigma - \Sigma)$ in terms of the effective rank $\ttr(\Sigma) = \Tr(\Sigma) / \|\Sigma\|$ under mild assumptions, completely eliminating dependence on the ambient dimension $d$. For example, \cite[Proposition A.3]{bunea15} prove that, for any $\delta \in (0, 1)$, with probability at least $1 - \delta$, we have
\begin{equation}
    \label{eq:dimension-free}
    \left\| \widehat\Sigma - \Sigma \right\|_{\Fr}
    \lesssim \Tr(\Sigma) \left( \sqrt{\frac{\log(2/\delta)}n} \vee \frac{\log(2/\delta)}n \right).
\end{equation}
Here and further in the paper, the sign $\lesssim$ stands for an inequality up to an absolute multiplicative constant. This brings us to a question, whether it is possible to obtain sharper dimension-free bounds when $\Sigma$ admits a Kronecker product expansion \eqref{eq:sigma_kronecker}. The following motivating example gives a positive answer in the simple case $K = 1$.

\begin{Ex}
    \label{ex:kronecker_rank_one}
    Let $\Sigma$ have a Kronecker rank one, that is, it admits a representation \eqref{eq:sigma_kronecker_simple} for some
    positive semidefinite matrices $\Phi \in \R^{p \times p}$ and $\Psi \in \R^{q \times q}$. Let $\X_1, \dots, \X_n \in \R^{q \times p}$ be the matrices obtained by reshaping\footnote{More precisely, it holds that $\vec(\X_i) = \bX_i$ for all $i \in \{1, \dots, n\}$, where $\vec$ is the vectorization operator, stacking the matrix columns together (see the definition in Appendix \ref{sec:kronecker}).} of the vectors $\bX_1, \dots, \bX_n$, respectively. Consider the estimates
    \[
        \widehat \Phi = \frac1n \sum\limits_{i = 1}^n \X_i^\top \X_i
        \quad \text{and} \quad
        \widehat \Psi = \frac1n \sum\limits_{i = 1}^n \X_i \X_i^\top.
    \]
    As shown by \cite{leng18}, $\widehat \Phi$ and $\widehat \Psi$ are unbiased estimates of $\Tr(\Psi) \Phi$ and $\Tr(\Phi) \Psi$, respectively. If, for any vectors $\bu \in \R^p$, $\bv \in \R^q$ and any $i \in \{1, \dots, n\}$, it holds that
    \[
        \frac1{\Tr(\Psi)} \left\| \bu^\top \X_i^\top \X_i \bu \right\|_{\psi_1} \lesssim \bu^\top \Phi \bu
        \quad \text{and} \quad
        \frac1{\Tr(\Phi)} \left\| \bv^\top \X_i \X_i^\top \bv \right\|_{\psi_1} \lesssim \bv^\top \Psi \bv,
    \]
    where $\|\cdot\|_{\psi_1}$ denotes the Orlicz norm (see the notation section below), then, similarly to the proof of Proposition A.3 from \citep{bunea15}, we can show that
    \[
        \left\|\widehat \Phi / \Tr(\Psi) - \Phi \right\|_{\Fr} \lesssim  \Tr(\Phi) \left( \sqrt{\frac{\log(2/\delta)}n} \vee \frac{\log(2/\delta)}n \right)
    \]
    and
    \[
        \left\|\widehat \Psi / \Tr(\Phi) - \Psi \right\|_{\Fr} \lesssim  \Tr(\Psi) \left( \sqrt{\frac{\log(2/\delta)}n} \vee \frac{\log(2/\delta)}n \right)  
    \]
    simultaneously on an event with probability at least $1 - \delta$. The trace product $\Tr(\Phi) \Tr(\Psi) = \Tr(\Sigma)$ can be easily estimated by $(\|\bX_1\|^2 + \dots + \|\bX_n\|^2) / n$. Then the estimate
    \begin{equation}
        \label{eq:leng-pan_estimate}
        \breve \Sigma = \frac{\widehat \Phi \otimes \widehat \Psi}{\frac1n \sum\limits_{i=1}^n \|\bX_i\|^2},
    \end{equation}
    suggested by \cite{leng18}, satisfies the following high probability bound:
    \begin{align*}
        \left\|\breve \Sigma - \Sigma \right\|_{\Fr}
        &
        \lesssim \left(\|\Psi\| \Tr(\Phi) + \|\Phi\| \Tr(\Psi) \right) \left( \sqrt{\frac{\log(2/\delta)}n} \vee \frac{\log(2/\delta)}n \right)
        \\&
        = \|\Sigma\| \left(\ttr(\Phi) + \ttr(\Psi) \right) \left( \sqrt{\frac{\log(2/\delta)}n} \vee \frac{\log(2/\delta)}n \right),
    \end{align*}
    where $\ttr(\Phi) = \Tr(\Phi) / \|\Phi\|$ and $\ttr(\Psi) = \Tr(\Psi) / \|\Psi\|$, while \eqref{eq:dimension-free} implies that
    \begin{align*}
        \left\|\widehat \Sigma - \Sigma \right\|_\Fr
        &
        \lesssim  \Tr(\Sigma) \left( \sqrt{\frac{\log(2/\delta)}n} \vee \frac{\log(2/\delta)}n \right)
        \\&
        = \|\Sigma\| \ttr(\Phi) \ttr(\Psi) \left( \sqrt{\frac{\log(2/\delta)}n} \vee \frac{\log(2/\delta)}n \right).
    \end{align*}
    Hence, the Kronecker structural assumption helped us to significantly improve the upper bound \eqref{eq:dimension-free} by replacing the product of effective ranks with their sum.
\end{Ex}

\bigskip

Unfortunately, the estimate $\breve \Sigma$ from Example \ref{ex:kronecker_rank_one} significantly exploits the product structure \eqref{eq:sigma_kronecker_simple} of $\Sigma$, which allows to estimate $\Phi$ and $\Psi$ separately. This approach has no direct extension on the case of larger Kronecker rank. For this reason, the case, when $\Sigma$ has a form \eqref{eq:sigma_kronecker} with $K > 1$, requires much more efforts.

\medskip

\noindent\textbf{Contribution.}\quad Our main contribution is a sharper upper bound on the performance of a regularized permuted least squares estimate $\widetilde \Sigma$ (defined in \eqref{eq:estimate}), which is widely used for covariance estimation in the presence of a Kronecker product structure. Under reasonable assumptions on the distribution of $\bX, \bX_1, \dots, \bX_n$, we show that, if $\Sigma$ has a form \eqref{eq:sigma_kronecker}, then
\[
    \left\| \widetilde \Sigma - \Sigma \right\|_\Fr^2
    \lesssim K \left( \sum\limits_{j = 1}^K \|\Phi_j\| \|\Psi_j\| \right)^2
    \cdot \left( \frac{1}{n} \max\limits_{1 \leq j \leq K} \ttr(\Phi_j)^2 + \frac{1}{n} \max\limits_{1 \leq j \leq K} \ttr(\Psi_j)^2 + \frac{\log(1 / \delta)}n \right)
\]
with probability at least $(1 - \delta)$. To our knowledge, this is the first dimension-free upper bound, capturing the Kronecker structure of the population covariance.

\medskip

\noindent\textbf{Notation.}\quad Throughout the paper, $\widehat \Sigma$ stands for the sample covariance. The identity matrix of size $(d \times d)$ and the zero matrix are denoted by $I_d$ and $O$, respectively. For any matrix $A$, $\|A\|$, $\|A\|_*$, and $\|A\|_\Fr$ stand for the operator, nuclear, and Frobenius norms of $A$, respectively. 
Given arbitrary matrices $A \in \R^{p \times q}$ and $B \in \R^{r \times s}$, $A \otimes B$ denotes their Kronecker product:
\[
    A \otimes B =
    \begin{pmatrix}
        a_{11} B & \dots & a_{1q} B \\
        \vdots & \ddots & \vdots \\
        a_{p1} B & \dots & a_{pq} B 
    \end{pmatrix}
    \in \mathbb{R}^{pr \times qs}.
\]
For any matrix $U \in \R^{k \times \ell}$ with columns $\bu_1, \dots, \bu_\ell$, $\vec(U) = (\bu_1^\top, \dots, \bu_\ell^\top)^\top$ is a vector obtained by stacking the columns of $U$ together. We also use the standard notation $\ttr(M) = \Tr(M) / \|M\|$ for the effective rank of a symmetric positive semidefinite matrix $M$. Instead of $\min\{a, b\}$ and $\max\{a, b\}$, sometimes we will use $a \land b$ and $a \lor b$, respectively. The notation $f \lesssim g$ is equivalent to $f = \cO(g)$. Finally, of any random variable $\xi$ and any $s \geq 1$, the Orlicz norm $\|\xi\|_{\psi_s}$ is defined as
\[
    \|\xi\|_{\psi_s} = \inf\left\{ t > 0 : \E e^{|\xi|^s / t^s} \leq 2 \right\}.
\]
The boldface font is reserved for vectors, while matrices and scalars are displayed in a regular font.

\medskip

\noindent\textbf{Paper structure.}\quad
The rest of the paper is organized as follows. In Section \ref{sec:main}, we present our main result, Theorem \ref{th:fr_norm_bound}, discuss its novelty and compare with related work. Sections \ref{sec:th_fr_norm_bound_proof} and \ref{sec:lem_op_norm_bound_proof} are devoted to the proof of Theorem \ref{th:fr_norm_bound} and its main ingredient, Lemma \ref{lem:op_norm_bound}. We also provide some useful information about Kronecker product and the vectorization operator $\vec$ in Appendix \ref{sec:kronecker}. The proofs of some technical results are also deferred to Appendix.

\section{Main results}
\label{sec:main}

This section contains the main result of this paper, but first, we have to introduce auxiliary definitions. One of the central objects when dealing with a covariance $\Sigma$ of the form \eqref{eq:sigma_kronecker} is the rearrangement operator $\ro$, defined below.

\begin{Def}[rearrangement operator]
    Let $M$ be a matrix of size $pq \times pq$ with a block structure
    \[
        M =
        \begin{pmatrix}
            M(1, 1) & \dots & M(1, p) \\
            \vdots & \ddots & \vdots \\
            M(p, 1) & \dots & M(p, p)
        \end{pmatrix},
        \quad
        \text{where $M(i, j) \in \R^{q \times q}$ for all $i, j \in \{1, \dots, p\}$.}
    \]
    The rearrangement operator $\ro : \R^{pq \times pq} \rightarrow \R^{p^2 \times q^2}$ is defined by the following identity: the $((j - 1)p + i)$-th row of $\ro(M)$ is $\vec(M(i, j))^\top$.
\end{Def}
The importance of the rearrangement operator becomes clear from the properties formulated in the next proposition. 
\begin{Prop}
    \label{prop:ro_properties}
    The following holds.
    \begin{itemize}
        \item $\ro$ is a linear operator.
        \item $\ro$ defines an isometry between the Euclidean spaces $(\R^{pq \times pq}, \|\cdot\|_\Fr)$ and $(\R^{p^2 \times q^2}, \|\cdot\|_\Fr)$. As a consequence, there exists an inverse map $\ro^{-1}$.
        \item If $M = A \otimes B$, where $A \in \R^{p \times p}$ and $B \in \R^{q \times q}$, then $\ro(M) = \vec(A) \vec(B)^\top$.
    \end{itemize}
\end{Prop}
The first and the second properties follow immediately from the definition of $\ro$. The third one is well-known in the literature (see, for instance, the proof of Theorem 1 in \citep{vanloan93}). The rearrangement operator reduces our problem to low-rank matrix estimation. Indeed, if the covariance matrix $\Sigma$ can be represented a sum of Kronecker products \eqref{eq:sigma_kronecker}, then, due to Proposition \ref{prop:ro_properties},
\[
    \ro(\Sigma)
    = \sum\limits_{j = 1}^K \ro(\Phi_j \otimes \Psi_j)
    = \sum\limits_{j = 1}^K \vec(\Phi_j) \vec(\Psi_j)^\top.
\]
In other words, if $\Sigma$ has a Kronecker rank $K$, then $\ro(\Sigma)$ has the same algebraic rank. This observation allows us to construct estimates based on the singular value decomposition of $\ro(\widehat\Sigma)$, as in \citep{tsiligkaridis13, masak22}.
For instance, \cite{masak22} considered an estimate based on $K$ principal components of $\ro(\widehat \Sigma)$:
\begin{equation}
    \label{eq:hard_thresholding}
    \widecheck \Sigma = \ro^{-1}( \widecheck R),
    \quad \text{where} \quad
    \widecheck R = \sum\limits_{j = 1}^{K} \widehat \sigma_j \widehat \bu_j \widehat \bv_j^\top,
\end{equation}
where the values $\widehat \sigma_1 \geq \dots \geq \widehat \sigma_K$ and the vectors $\widehat \bu_1, \widehat \bv_1, \dots, \widehat \bu_K, \widehat \bv_K$ are defined from the singular value decomposition
\[
    \ro(\widehat\Sigma) = \sum\limits_{j = 1}^{p^2 \land q^2} \widehat \sigma_j \widehat \bu_j \widehat \bv_j^\top.
\]
Let us note that the estimate \eqref{eq:hard_thresholding} solves the optimization problem
\[
    \widecheck R \in \argmin\limits_{R \in \R^{p^2 \times q^2}} \left\|R - \ro(\widehat \Sigma) \right\|_\Fr^2 + \lambda \, \rank(R),
    \quad \text{where $\widehat\sigma_{K + 1}^2 \leq \lambda \leq \widehat\sigma_K^2$.}
\]
In the present paper, we study the permuted least squares estimate suggested in \citep{tsiligkaridis13}, where the authors replaced $\rank(R)$ by nuclear norm penalization:
\begin{equation}
    \label{eq:estimate}
    \widetilde \Sigma = \ro^{-1}(\widetilde R),
    \quad \text{where} \quad
    \widetilde R \in \argmin\limits_{R \in \R^{p^2 \times q^2}} \left\{ \left\|R - \ro(\widehat \Sigma) \right\|_\Fr^2 + \lambda \|R\|_* \right\}.
\end{equation}
Similarly to $\widecheck R$, the estimate $\widetilde R$ admits an explicit representation given by soft thresholding:
\[
    \widetilde R = \sum\limits_{j = 1}^{p^2 \land q^2} \big( \widehat \sigma_j - \lambda / 2 \big)_+ \widehat \bu_j \widehat \bv_j^\top.
\]
Our goal is to establish a dimension-free upper bound on the Frobenuis norm of $(\widetilde \Sigma - \Sigma)$. Usually, results of such kind require some additional properties of the distribution of $\bX, \bX_1, \dots, \bX_n$. We impose the following assumption.

\begin{As}
    \label{as:orlicz}
    There exists $\omega > 0$, such that the standardized random vector $\Sigma^{-1/2} \bX$ satisfies the inequality
    \begin{equation}
        \label{eq:orlicz_norm_inequality}
        \log \E \exp \left\{ (\Sigma^{-1/2} \bX)^\top V (\Sigma^{-1/2} \bX) - \Tr(V) \right\} \leq \omega^2 \|V\|_{\Fr}^2
    \end{equation}
    for all $V \in \R^{d \times d}$, such that $\|V\|_{\Fr} \leq 1 /  \omega$.
\end{As}
An equivalent condition appeared in the paper \citep{puchkin23}, where the authors assumed that
\begin{equation}
    \label{eq:orlicz_norm_equivalent}
    \left\|(\Sigma^{-1/2} \bX)^\top V (\Sigma^{-1/2} \bX) - \Tr(V) \right\|_{\psi_1} \lesssim \|V\|_{\Fr}
    \quad \text{for all $V \in \R^{d \times d}$.}
\end{equation}
The fact that \eqref{eq:orlicz_norm_equivalent} yields \eqref{eq:orlicz_norm_inequality} follows from, for instance, \citep[Lemma 2]{zhivotovskiy24}. In \citep{puchkin23}, the authors argued that 
Assumption \ref{as:orlicz} holds for a large class of distributions. In particular, they discussed that it is fulfilled for all random vectors, satisfying the Hanson-Wright inequality. The examples include the Gaussian distribution $\N(\bzero, \Sigma)$ and, more generally, all distributions, such that $\Sigma^{-1/2} \bX$ has independent components with finite $\psi_2$-norms. Besides, random vectors satisfying log-Sobolev inequality or having the convex concentration property \citep{adamczak15} also meet Assumption \ref{as:orlicz}. Since a covariance matrix with a Kronecker product structure \eqref{eq:sigma_kronecker} often appears in the context of matrix models, we would like to give relevant examples from this field. For instance, if we have a matrix model
\[
    \X = B \Y A^\top,
\]
where $B \in \R^{q \times q}$, $A \in \R^{p \times p}$, and $\Y \in \R^{q \times p}$ is a random matrix with i.i.d. standard Gaussian entries, then $\vec(\X) \sim \N \big(\bzero, (AA^\top) \otimes (BB^\top) \big)$ and, consequently, it fulfils Assumption \ref{as:orlicz}. In the proposition below, we provide another, more general example when the assumption is satisfied.

\begin{Prop}
    \label{prop:orlicz}
    Let
    \begin{equation}
        \label{eq:matrix_model}
        \X = \sum\limits_{j = 1}^K B_j \Y_j A_j^\top,
    \end{equation}
    where $\Y_1, \dots, \Y_K \in \R^{q \times p}$ are independent centered random matrices. Assume that there exist positive numbers $\omega_1, \dots, \omega_K$, such that for any $j \in \{1, \dots, K\}$
    \begin{equation}
        \label{eq:y_orlicz_norm}
        \log \E \exp \left\{ \vec(\Y_j)^\top V \vec(\Y_j) - \Tr(V) \right\} \leq \omega_j^2 \|V\|_{\Fr}^2
    \end{equation}
    for all $V \in \R^{d^2 \times d^2}$, satisfying the inequality $\|V\|_\Fr \leq 1 / \omega_j$. Then $\vec(\X)$ fulfils Assumption \ref{as:orlicz} with
    \[
        \Sigma = \sum\limits_{j = 1}^K (A_j A_j^\top) \otimes (B_j B_j^\top)
        \quad \text{and} \quad
        \omega \leq \left(2 \omega_{\max} \vee 16 \sqrt{\frac{2 + 2\omega_{\max}}{\log 2}} \right),
    \]
    where $\omega_{\max} = \max\limits_{1 \leq j \leq K} \omega_j$.
\end{Prop}
The proof of Proposition \ref{prop:orlicz} is quite technical. For this reason, we postpone it to Appendix \ref{sec:prop_orlicz_proof} and move to the main result.
\begin{Th}
    \label{th:fr_norm_bound}
    Assume that the centered random vector $\bX$ meets Assumption \eqref{as:orlicz} and that its covariance $\Sigma$ has a form \eqref{eq:sigma_kronecker}. Let us fix any $\delta \in (0, 1)$, satisfying the constraint
    \[
        \frac1{2n} \left( \max\limits_{1 \leq j \leq K} \ttr(\Phi_j)^2 + \max\limits_{1 \leq j \leq K} \ttr(\Psi_j)^2 \right)
        + \frac{\log(4 / \delta)}n
        \leq 1,
    \]
    and take any $\lambda \in \R$, such that
    \[
        \lambda
        \geq 2 \omega
        \left( \sum\limits_{j = 1}^K \|\Phi_j\| \|\Psi_j\| \right)
        \sqrt{\frac{13}{2n} \left( \max\limits_{1 \leq j \leq K} \ttr(\Phi_j)^2 + \max\limits_{1 \leq j \leq K} \ttr(\Psi_j)^2 \right) + \frac{13 \log(4 / \delta)}n}.
    \]
    Then, with probability at least $(1 - \delta)$, the permuted nuclear-norm penalized least squares estimate $\widetilde \Sigma$, defined in \eqref{eq:estimate}, fulfils the inequality
    \[
        \left\| \widetilde \Sigma - \Sigma \right\|_\Fr^2
        < \frac{3 \lambda^2 K}2.
    \]
\end{Th}

\begin{Rem}
    Theorem \ref{th:fr_norm_bound} requires a proper choice of $\lambda > 0$. To tune the parameter, we can fix a finite family $\Lambda = \{\lambda_0 \cdot 2^{-m} : 0 \leq m \leq M\}$, where $\lambda_0 > 0$ and $M \in \mathbb N$ are some constants, and perform model selection, using cross-validation or the random split procedure described in \cite[Section 3]{bickel08b}. If the chosen $\lambda$ meets the inequality
    \[
        \lambda
        \leq C \omega \left( \sum\limits_{j = 1}^K \|\Phi_j\| \|\Psi_j\| \right) \sqrt{\frac{13}{2n} \left( \max\limits_{1 \leq j \leq K} \ttr(\Phi_j)^2 + \max\limits_{1 \leq j \leq K} \ttr(\Psi_j)^2 \right) + \frac{13 \log(4 / \delta)}n}
    \]
    with an absolute constant $C > 2$, then, under the conditions of Theorem \ref{th:fr_norm_bound}, it holds that
    \[
        \left\| \widetilde \Sigma - \Sigma \right\|_\Fr^2
        \lesssim \omega^2 K
        \left( \sum\limits_{j = 1}^K \|\Phi_j\| \|\Psi_j\| \right)^2
        \left( \frac1n \max\limits_{1 \leq j \leq K} \ttr(\Phi_j)^2 + \frac1n \max\limits_{1 \leq j \leq K} \ttr(\Psi_j)^2 + \frac{\log(1 / \delta)}n \right)
    \]
    with probability at least $(1 - \delta)$. 
\end{Rem}

To our knowledge, Theorem \ref{th:fr_norm_bound} provides the first dimension-free high-probability upper bound on the accuracy of covariance estimation with a Kronecker product structure. The most relevant papers to ours are \citep{masak22} and \citep{tsiligkaridis13}, where the authors considered the estimates $\widecheck \Sigma$ and $\widetilde \Sigma$ (see \eqref{eq:hard_thresholding}, \eqref{eq:estimate}), based on the hard and soft thresholding, respectively. In \citep{masak22}, the authors assumed only the existence of the fourth moment of $\|\bX\|$ and proper decay of the singular values of $\ro(\Sigma)$. Under these mild conditions, they 
proved the consistency of $\widecheck \Sigma$. Unfortunately, the established rate of convergence includes the Frobenius norm of $(\widehat \Sigma - \Sigma)$, which is of order $\ttr(\Sigma) / \sqrt{n}$ in the worst-case scenario and, in contrast to Theorem \ref{th:fr_norm_bound}, shows no improvements compared to the unstructured case. In \citep{tsiligkaridis13}, the authors derived a non-asymptotic bound on the performance of the permuted least-squares estimate $\widetilde \Sigma$. Assuming the observations $\bX_1, \dots, \bX_n$ to have a Gaussian distribution $\N(\bzero, \Sigma)$, they showed that
\[
    \left\|\widetilde \Sigma - \Sigma \right\|_\Fr^2
    \lesssim \frac{K (p^2 + q^2 + \log M)}n \vee \frac{K (p^2 + q^2 + \log M)^2}{n^2}
\]
with probability at least $1 - M^{-c}$, where $c > 0$ is an absolute constant. In our paper, we derive a sharper dimension-free bound for a broader class of distributions, which remains valid in the extremely high-dimensional setup.

Let us elaborate on implications of Theorem \ref{th:fr_norm_bound} in the simpler case $K = 1$, which also got a considerable attention in the literature. In \citep{werner08}, the authors focused on the Gaussian case and showed the consistency of the maximum likelihood estimate and of non-iterative flip-flop algorithm. Later, \cite{zhou14} suggested the Gemini method and obtained high-probability dimension-dependent bounds on the accuracy of covariance estimation, improving over the results of \cite*{werner08}. The Gemini estimator also exploits the fact that $\bX_1, \dots, \bX_n$ have the Gaussian distribution. In contrast, \cite*{leng18} considered a much broader class of distributions with finite $48$-th moment. Under some technical assumptions, they derived a high-probability upper bound on the accuracy of covariance estimation in terms of both operator and Frobenius norms. Unfortunately, the analysis of \cite{zhou14} and \cite{leng18} imposes a severe restriction that the lowest eigenvalues of $\Phi$ and $\Psi$ from the decomposition \eqref{eq:sigma_kronecker_simple} should be bounded away from zero. This makes their theoretical guarantees vacuous in the high-dimensional setting. On the other hand, if the Kronecker rank of $\Sigma$ is equal to one, then Theorem \ref{th:fr_norm_bound} yields that
\[
    \left\| \widetilde \Sigma - \Sigma \right\|_\Fr^2
    \lesssim \omega^2 \|\Sigma\|^2
    \cdot \frac{\ttr(\Phi)^2 + \ttr(\Psi)^2 + \log(1 / \delta)}n
\]
with probability at least $(1 - \delta)$. Note that this bound has a slightly better dependence on $\log(1/\delta)$, than the rate of convergence from Example \ref{ex:kronecker_rank_one}.

\section{Proof of Theorem \ref{th:fr_norm_bound}}
\label{sec:th_fr_norm_bound_proof}

The first step in the proof of the rate of convergence of the permuted least squares estimate $\widetilde \Sigma$ is identical to the one of \cite{tsiligkaridis13}. With a proper choice of $\lambda$, we have the following oracle inequality.

\begin{Th}[\citet*{tsiligkaridis13}, Theorem 2]
    If $\lambda \geq 2 \|\ro(\widehat \Sigma - \Sigma) \|$, then
    \begin{equation}
        \label{eq:oracle_inequality}
        \left\| \widetilde R - \ro(\Sigma) \right\|_\Fr^2
        \leq \inf\limits_{R} \left\{ \left\| R - \ro(\Sigma) \right\|_\Fr^2 + \frac{(1 + \sqrt 2)^2}4 \lambda^2 \rank(R) \right\}.
    \end{equation}
\end{Th}
Obviously, since $\Sigma$ admits a representation \eqref{eq:sigma_kronecker}, the rank of the matrix
\begin{equation}
    \label{eq:ro_sigma_rank}
    \ro(\Sigma)
    = \ro\left( \sum\limits_{j = 1}^K \Phi_j \otimes \Psi_j \right)
    = \sum\limits_{j = 1}^K \ro\left( \Phi_j \otimes \Psi_j \right)
    = \sum\limits_{j = 1}^K \vec(\Phi_j) \vec(\Psi_j)^\top
\end{equation}
is equal to $K$. In view of \eqref{eq:ro_sigma_rank}, the oracle inequality \eqref{eq:oracle_inequality} can be rewritten in the form
\[
    \left\| \widetilde R - \ro(\Sigma) \right\|_\Fr^2
    \leq \frac{(1 + \sqrt 2)^2 \lambda^2 K}4,
    \quad 
    \text{provided that $\lambda \geq 2 \|\ro(\widehat \Sigma - \Sigma)\|$.}
\]
Our main technical result, leading to superior guarantees on the performance of $\widetilde \Sigma$ is the following dimension-free bound on the operator norm of $\ro(\widehat \Sigma - \Sigma)$.
\begin{Lem}
    \label{lem:op_norm_bound}
    Let Assumption \ref{as:orlicz} hold and let us fix any $\delta \in (0, 1)$, satisfying the inequality
    \begin{equation}
        \label{eq:delta_condition}
        \frac1{2n} \left( \max\limits_{1 \leq j \leq K} \ttr(\Phi_j)^2 + \max\limits_{1 \leq j \leq K} \ttr(\Psi_j)^2 \right)
        + \frac{\log(4 / \delta)}n
        \leq 1.    
    \end{equation}
    Then, with probability at least $(1 - \delta)$, it holds that
    \[
        \left\| \ro(\widehat \Sigma - \Sigma) \right\|
        \leq \omega
        \left( \sum\limits_{j = 1}^K \|\Phi_j\| \|\Psi_j\| \right)
        \sqrt{\frac{13}{2n} \left( \max\limits_{1 \leq j \leq K} \ttr(\Phi_j)^2 + \max\limits_{1 \leq j \leq K} \ttr(\Psi_j)^2 \right) + \frac{13 \log(4 / \delta)}n}.
    \]
\end{Lem}
We present the proof of Lemma \ref{lem:op_norm_bound} in Section \ref{sec:lem_op_norm_bound_proof} below. With Lemma \ref{lem:op_norm_bound}, we easily finish the proof of Theorem \ref{th:fr_norm_bound}. Let $\lambda$ be any number, satisfying the inequality
\[
    \lambda
    \geq 2 \omega
    \left( \sum\limits_{j = 1}^K \|\Phi_j\| \|\Psi_j\| \right)
    \sqrt{\frac{13}{2n} \left( \max\limits_{1 \leq j \leq K} \ttr(\Phi_j)^2 + \max\limits_{1 \leq j \leq K} \ttr(\Psi_j)^2 \right) + \frac{13 \log(4 / \delta)}n}.
\]
Then, according to Lemma \ref{lem:op_norm_bound}, $\lambda \geq 2 \| \ro(\widehat \Sigma - \Sigma) \|$ on an event of probability at least $(1 - \delta)$. Hence, on this event, it holds that
\[
    \left\| \widetilde \Sigma - \Sigma \right\|_\Fr^2
    = \left\| \widetilde R - \ro(\Sigma) \right\|_\Fr^2
    \leq \frac{(1 + \sqrt 2)^2 \lambda^2 K}4
    = \frac{(3 + 2\sqrt 2) \lambda^2 K}4
    < \frac{3 \lambda^2 K}2.
\]
\endproof

\section{Proof of Lemma \ref{lem:op_norm_bound}}
\label{sec:lem_op_norm_bound_proof}

The present section is devoted to the proof of our main technical result. For convenience, we split it into several steps, deferring auxiliary derivations to Appendix.

\medskip

\noindent
\textbf{Step 1: reduction to an empirical process.}
\quad
We start with representing the operator norm of $\ro(\widehat \Sigma - \Sigma)$ in an appropriate form. 
\begin{Lem}
    \label{lem:ro_representation}
    Let $\X, \X_1, \dots, \X_n$ be random matrices of size $(q \times p)$, such that $\vec(\X) = \bX$ and $\vec(\X_i) = \bX_i$ for all $i \in \{1, \dots, n\}$. Then it holds that
    \begin{equation}
        \label{eq:ro_representation}
        \left\| \ro(\widehat \Sigma - \Sigma) \right\|
        = \sup\limits_{\substack{U \in \R^{p \times p}, V \in \R^{q \times q}\\ \|U\|_\Fr = \|V\|_\Fr = 1}} \left\{ \frac1n \sum\limits_{i = 1}^n \Tr\left( \X_i^\top V^\top \X_i U \right) - \E \Tr\left( \X^\top V^\top \X U \right) \right\}.
    \end{equation}
\end{Lem}
The proof of Lemma \ref{lem:ro_representation} is postponed to Appendix \ref{sec:lem_ro_representation_proof}. Now we are in position to use powerful tools from the empirical process theory to derive a dimension-free bound on the supremum in the right-hand side of \eqref{eq:ro_representation}. The main ingredient of our proof is the following PAC-Bayesian variational inequality (see, for instance, \citep[Proposition 2.1]{catoni17}). 

\begin{Lem}
    \label{lem:pac-bayes}
    Let $ \bX, \bX_1, \dots, \bX_n$ be i.i.d. random elements on a measurable space $\cX$. Let $\Theta$ be a parameter space equipped with a measure $\mu$ (which is also referred to as prior). Let $f : \cX \times \Theta \rightarrow \R$. Then, with probability at least $1 - \delta$, it holds that
    \[
        \E_{\btheta \sim \rho} \frac1n \sum\limits_{i = 1}^n f(\bX_i, \btheta)
        \leq \E_{\btheta \sim \rho} \log \E_\bX e^{f(\bX, \btheta)} + \frac{\KL(\rho, \mu) + \log(1 / \delta)}n
    \]
    simultaneously for all $\rho \ll \mu$.
\end{Lem}
Recently, \cite{zhivotovskiy24} used a similar approach to prove that the operator norm of $(\widehat \Sigma - \Sigma)$ is of order $\cO(\sqrt{\ttr(\Sigma) / n})$ with high probability. Unfortunately, the rate of convergence $\cO(\sqrt{\ttr(\Sigma) / n})$ is suboptimal in our case, because of the additional structure of the population covariance. In the rest of the proof, we deduce a sharper high-probability upper bound
\[
    \|\ro(\widehat \Sigma - \Sigma)\|
    \lesssim \sqrt{\frac1n \left(\max\limits_{1 \leq j \leq K} \ttr(\Phi_j)^2 + \max\limits_{1 \leq j \leq K} \ttr(\Psi_j)^2 \right)},
\]
leading to better guarantees on the accuracy of estimation of $\Sigma$.

\medskip

\noindent\textbf{Step 2: a variational inequality.}\quad On this step, we specify the prior and posterior distributions in order to apply Lemma \ref{lem:pac-bayes} to our problem. In our case, both $\mu$ and $\rho$ will be absolutely continuous with respect to the Lebesgue measure, so we identify them with the corresponding densities. Let $\mu$ be the density of a Gaussian measure on $\R^{p \times p} \times \R^{q \times q}$:
\begin{equation}
    \label{eq:prior_density}
    \mu(X, Y) = \frac{\alpha^{p^2 / 2} \beta^{q^2 / 2}}{(2 \pi)^{p^2 / 2 + q^2 / 2}} \exp\left\{ -\frac{\alpha}2 \|X\|_{\Fr}^2 - \frac{\beta}2 \|Y\|_{\Fr}^2 \right\},
    \quad \text{where $X \in \R^{p \times p}$, $Y \in \R^{q \times q}$,}
\end{equation}
and
\begin{equation}
    \label{eq:alpha_beta}
    \alpha = \max\limits_{1 \leq j \leq K} \ttr(\Phi_j)^2,
    \quad
    \beta = \max\limits_{1 \leq j \leq K} \ttr(\Psi_j)^2.
\end{equation}
For any $U \in \R^{p \times p}$ and $V \in \R^{q \times q}$, define a set
\begin{align}
    \label{eq:support}
    \Upsilon(U, V)
    &\notag
    = \Bigg\{ (X, Y) \in \R^{p \times p} \times \R^{q \times q} :
    \left\| \Sigma^{1/2} \big( (X - U) \otimes (Y - V) \big) \Sigma^{1/2} \right\|_\Fr^2 \leq \frac{4 \Tr(\Sigma)^2}{\alpha \beta},
    \\&\qquad
    \left\| \left( \sum\limits_{j = 1}^K \|\Psi_j\| \Phi_j \right)^{1/2} (X - U) \left( \sum\limits_{j = 1}^K \|\Psi_j\| \Phi_j \right)^{1/2} \right\|_{\Fr}^2 \leq \frac{4}{\alpha} \left( \sum\limits_{j = 1}^K \|\Psi_j\| \Tr(\Phi_j) \right)^{2},
    \\&\qquad\notag
    \left\| \left( \sum\limits_{j = 1}^K \|\Phi_j\| \Psi_j \right)^{1/2} (Y - V) \left( \sum\limits_{j = 1}^K \|\Phi_j\| \Psi_j \right)^{1/2} \right\|_{\Fr}^2 \leq \frac{4}{\beta} \left( \sum\limits_{j = 1}^K \|\Phi_j\| \Tr(\Psi_j) \right)^{2}
    \Bigg\}
\end{align}
and a corresponding posterior density $\rho_{U, V}$, supported on $\Upsilon(U, V)$:
\begin{equation}
    \label{eq:posterior_density}
    \rho_{U, V}(X, Y) =
    \begin{cases}
        \mu(X - U, Y - V) / \ttZ, \quad \text{if $(X, Y) \in \Upsilon(U, V)$},\\
        0, \quad \text{otherwise,}
    \end{cases}
\end{equation}
where $\ttZ$ is a normalizing constant. It is not hard to show that $\ttZ$ is at least $1/4$.
\begin{Lem}
    \label{lem:z}
    With the notations introduced above, it holds that $\ttZ \geq 1/4$.
\end{Lem}
The proof of Lemma \ref{lem:z} is deferred to Appendix \ref{sec:lem_z_proof}. Let $\nu > 0$ be a positive constant to be defined later. Considering a pair $(U, V) \in \R^{p \times p} \times \R^{q \times q}$ as an element of the parameter space
\[
    \Theta = \left\{ (U, V) \in \R^{p \times p} \times \R^{q \times q} \right\},
\]
we conclude that, according to Lemma \ref{lem:pac-bayes},
\begin{align}
    \label{eq:pac-bayes}
    &\notag
    \sup\limits_{\substack{U \in \R^{p \times p}, V \in \R^{q \times q}\\ \|U\|_\Fr = \|V\|_\Fr = 1}} \left\{ \frac{\nu}n \sum\limits_{i = 1}^n \Tr\left( \X_i^\top V^\top \X_i U \right) - \nu \E_\X \Tr\left( \X^\top V^\top \X U \right) \right\}
    \\&
    = \sup\limits_{\substack{U \in \R^{p \times p}, V \in \R^{q \times q}\\ \|U\|_\Fr = \|V\|_\Fr = 1}} \left\{ \frac{\nu}n \sum\limits_{i = 1}^n \E_{(P, Q) \sim \rho_{U, V}} \Tr\left( \X_i^\top Q^\top \X_i P \right) - \nu \E_{(P, Q) \sim \rho_{U, V}} \E_\X \Tr\left( \X^\top Q^\top \X P \right) \right\}
    \\&\notag
    \leq \sup\limits_{\substack{U \in \R^{p \times p}, V \in \R^{q \times q}\\ \|U\|_\Fr = \|V\|_\Fr = 1}} \Bigg\{ \E_{(P, Q) \sim \rho_{U, V}} \log \left[ \E_\X \exp\left\{ \nu \Tr\left( \X^\top Q^\top \X P \right) - \nu \E_\X \Tr\left( \X^\top Q^\top \X P \right) \right\} \right]
    \\&\hspace{4in}\notag
    + \frac{\KL(\rho_{U, V}, \mu) + \log(1 / \delta)}n \Bigg\}.
\end{align}
with probability at least $(1 - \delta)$.

\medskip

\noindent
\textbf{Step 3: a bound on the exponential moment.}
\quad
Our next goal is to bound the exponential moment
\[
    \E_\X \exp\left\{ \nu \Tr\left( \X^\top Q^\top \X P \right) - \nu \E_\X \Tr\left( \X^\top Q^\top \X P \right) \right\}.
\]
We begin with the following auxiliary result.
\begin{Lem}
    \label{lem:support}
    Let us fix any $U \in \R^{p \times p}$ and $V \in \R^{q \times q}$, such that $\|U\|_\Fr = \|V\|_\Fr = 1$. Then
    \[
        \left\|\Sigma^{1/2} (P \otimes Q) \Sigma^{1/2} \right\|_{\Fr}^2
        \leq 13 \left( \sum\limits_{j = 1}^K \|\Phi_j\| \|\Psi_j\| \right)^2
        \quad \text{$\rho_{U,V}$-almost surely.}
    \]
\end{Lem}
The proof of Lemma \ref{lem:support} is moved to Appendix \ref{sec:lem_support_proof}.
Applying this lemma and using the fact that
$\bX = \vec(\X)$ satisfies Assumption \ref{as:orlicz}, we obtain that
\begin{align}
    \label{eq:exp_moment_bound}
    &\notag
    \log \E_\X \exp\left\{ \nu \Tr\left( \X^\top Q^\top \X P \right) - \nu \E_\X \Tr\left( \X^\top Q^\top \X P \right) \right\}
    \\&\notag
    = \log \E_\X \exp\left\{ \nu \vec(\X)^\top (P \otimes Q) \vec(\X) - \nu \E_\X \vec(\X)^\top (P \otimes Q) \vec(\X) \right\}
    \\&
    = \log \E_{\bX} \exp\left\{ \nu \bX^\top (P \otimes Q) \bX - \nu \E_{\bX} \bX^\top (P \otimes Q) \bX \right\}
    \\&\notag
    \leq \omega^2 \nu^2 \left\|\Sigma^{1/2} (P \otimes Q) \Sigma^{1/2} \right\|_{\Fr}^2
    \\&\notag
    \leq 13 \omega^2 \nu^2 \left( \sum\limits_{j = 1}^K \|\Phi_j\| \|\Psi_j\| \right)^2
\end{align}
for all $\nu$, such that
\begin{equation}
    \label{eq:nu_condition}
    13 \omega^2 \nu^2 \left( \sum\limits_{j = 1}^K \|\Phi_j\| \|\Psi_j\| \right)^2 \leq 1,
\end{equation}
and for all $U \in \R^{p \times p}$, $V \in \R^{q \times q}$ with unit Frobenius norms.

\medskip

\noindent
\textbf{Step 4: a bound on the Kullback-Leibler divergence.} 
\quad
We proceed with an upper bound on $\KL(\rho_{U, V}, \mu)$, where the densities $\mu$ and $\rho_{U, V}$ are defined in \eqref{eq:prior_density} and \eqref{eq:posterior_density}, respectively. It holds that
\[
    \KL(\rho_{U, V}, \mu)
    = \int\limits_{\Upsilon(U, V)} \log\frac{\rho_{U, V}(X, Y)}{\mu(X, Y)} \rho_{U, V}(X, Y) \dd X \dd Y.
\]
Substituting $X$ and $Y$ by $X + U$ and $Y + V$, respectively, we immediately obtain that
\begin{align*}
    \KL(\rho_{U, V}, \mu)
    &
    = \int\limits_{\Upsilon(O, O)} \log\frac{\rho_{U, V}(X + U, Y + V)}{\mu(X + U, Y + V)} \rho_{U, V}(X + U, Y + V) \dd X \dd Y
    \\&
    = \int\limits_{\Upsilon(O, O)} \log\frac{\rho_{O, O}(X, Y)}{\mu(X + U, Y + V)} \rho_{O, O}(X, Y) \dd X \dd Y
    \\&
    = \log\frac1{\ttZ}
    + \int\limits_{\Upsilon(O, O)} \left( \alpha \Tr(U^\top X) + \beta \Tr(V^\top Y) \right)
    \rho_{O, O}(X, Y) \dd X \dd Y
    \\&\quad
    + \frac12 \int\limits_{\Upsilon(O, O)} \left( \alpha \|U\|_\Fr^2 + \beta \|V\|_\Fr^2 \right)
    \rho_{O, O}(X, Y) \dd X \dd Y
    \\&
    = \log\frac1{\ttZ} + \frac{\alpha \|U\|_\Fr^2 + \beta \|V\|_\Fr^2}2.
\end{align*}
Here we took into account that $\Upsilon(O, O)$ and $\rho_{O, O}$ are symmetric around zero and, hence,
\[
    \int\limits_{\Upsilon(O, O)} \Tr(U^\top X)
    \rho_{O, O}(X, Y) \dd X \dd Y
    = 0,
    \quad
    \int\limits_{\Upsilon(O, O)} \Tr(V^\top Y)
    \rho_{O, O}(X, Y) \dd X \dd Y
    = 0.
\]
Applying Lemma \ref{lem:z}, we conclude that
\begin{equation}
    \label{eq:kl_bound}
    \KL(\rho_{U, V}, \mu)
    \leq 2 \log 2 + \frac{\alpha + \beta}2
    = 2\log2 + \frac12 \left( \max\limits_{1 \leq j \leq K} \ttr(\Phi_j)^2 + \max\limits_{1 \leq j \leq K} \ttr(\Psi_j)^2 \right)
\end{equation}
for all $U \in \R^{p \times p}$, $V \in \R^{q \times q}$, such that $\|U\|_\Fr = \|V\|_\Fr = 1$.

\medskip

\noindent
\textbf{Step 4: final bound.}
\quad
Summing up \eqref{eq:ro_representation}, \eqref{eq:pac-bayes}, \eqref{eq:exp_moment_bound}, and \eqref{eq:kl_bound}, we obtain that
\begin{align*}
    \nu \left\| \ro(\widehat \Sigma - \Sigma) \right\|
    &
    \leq 13 \omega^2 \nu^2 \left( \sum\limits_{j = 1}^K \|\Phi_j\| \|\Psi_j\| \right)^2 
    \\&\quad
    + \frac1{2n} \left( \max\limits_{1 \leq j \leq K} \ttr(\Phi_j)^2 + \max\limits_{1 \leq j \leq K} \ttr(\Psi_j)^2 \right)
    + \frac{\log(4 / \delta)}n
\end{align*}
with probability at least $(1 - \delta)$. Let us take
\[
    \nu = \frac1{\omega \sqrt{13}} \left( \sum\limits_{j = 1}^K \|\Phi_j\| \|\Psi_j\| \right)^{-1} \sqrt{\frac1{2n} \left( \max\limits_{1 \leq j \leq K} \ttr(\Phi_j)^2 + \max\limits_{1 \leq j \leq K} \ttr(\Psi_j)^2 \right) + \frac{\log(4 / \delta)}n}.
\]
Note that the requirement \eqref{eq:delta_condition} ensures that such choice of $\nu$ fulfils \eqref{eq:nu_condition}. Hence, with probability at least $(1 - \delta)$, it holds that
\[
    \left\| \ro(\widehat \Sigma - \Sigma) \right\|
    \leq \omega \left( \sum\limits_{j = 1}^K \|\Phi_j\| \|\Psi_j\| \right) \sqrt{\frac{13}{2n} \left( \max\limits_{1 \leq j \leq K} \ttr(\Phi_j)^2 + \max\limits_{1 \leq j \leq K} \ttr(\Psi_j)^2 \right) + \frac{13 \log(4 / \delta)}n}.
\]
\endproof

\section{Conclusion and open problems}

We showed that, if the covariance matrix $\Sigma$ can be represented in the form \eqref{eq:sigma_kronecker}, then it is possible to obtain a non-asymptotic dimension-free upper bound on the Frobenius distance between $\Sigma$ and the penalized permuted least squares estimate $\widetilde \Sigma$ under reasonable conditions. The rate of convergence is sharper compared to the standard unstructured case. There is an open question whether similar improvements are possible if we measure the estimation performance in terms of the operator norm. In particular, we are not aware of any (even dimension-dependent) bounds on $\|\widetilde \Sigma - \Sigma\|$. We also leave open the question if the rate of convergence in Theorem \ref{th:fr_norm_bound} is optimal in the minimax sense. Another interesting question is how the result of Theorem \ref{th:fr_norm_bound} will change if we consider a misspecified model
\[
    \Sigma
    = \Phi_1 \otimes \Psi_1 + \ldots + \Phi_K \otimes \Psi_K + E,
\]
where $E \in \R^{d \times d}$ is a (possibly unstructured) remainder. In \citep{greenewald15}, the authors studied the case when $E$ is a sparse matrix and suggested a robust Kronecker PCA procedure for estimation of $\Sigma$ in this situation. Finally, one can consider more complex models with higher-order Kronecker products \citep{pouryazdian16, hafner20, mccormack23}.

\acks{The article was prepared within the framework of the HSE University Basic Research Program.}

\bibliography{references.bib}

\appendix

\section{Auxiliary results: Kronecker product and its properties}
\label{sec:kronecker}

The goal of this section is to introduce key properties of Kronecker products that we utilize in this paper. For more details see, for instance, \citep{golub2013matrix}. 
We start with a basic fact also referred to as the mixed-product property:
\begin{equation}
    \label{eq:mixed_product}
    (A \otimes B) (C \otimes D) = (AC) \otimes (BD),
\end{equation}
The equality \eqref{eq:mixed_product} has plenty of implications. In particular, if $A = U^\top \Lambda U$ and $B = V^\top M V$ are eigendecompositions of symmetric matrices $A \in \R^{p \times p}$ and $B \in \R^{q \times q}$, then
\[
    A \otimes B
    = (U^\top \Lambda U) \otimes (V^\top M V)
    = (U \otimes V)^\top (\Lambda \otimes M) (U \otimes V)
\]
is the eigendecomposition of their Kronecker product. For this reason, we have
\begin{equation}
    \label{eq:kronecker_eigenvalues}
    \Tr(A \otimes B) = \Tr(A) \Tr(B),
    \quad
    \|A \otimes B\| = \|A\| \|B\|,
    \quad
    \text{and}
    \quad
    \|A \otimes B\|_{\Fr} = \|A\|_{\Fr} \|B\|_{\Fr}.
\end{equation}

Further properties of the Kronecker product, presented in this section, are related to the vectorization operator $\vec$ (see our notation in Section \ref{sec:introduction}). Let us recall that, for any matrix $U$, $\vec(U)$ is a vector obtained by stacking the columns of $U$ together. In our proofs, we will extensively use the following identities, which hold whenever the matrix products in the left- and the right-hand side are well defined: 
\begin{equation}
    \label{eq:kronecker_vec}
    (A \otimes B) \vec(U) = \vec(B U A^\top),
\end{equation}
\begin{equation}
    \label{eq:kronecker_trace}
    \Tr(V^\top B U A^\top) = \vec(V)^\top \vec(B U A^\top) = \vec(V)^\top (A \otimes B) \vec(U).
\end{equation}

\section{Proof of Proposition \ref{prop:orlicz}}
\label{sec:prop_orlicz_proof}

Let us fix any $V \in \R^{d^2 \times d^2}$ and consider the exponential moment
\[
    \E \exp\left\{ \vec(\X)^\top V \vec(\X) - \Tr(V \Sigma) \right\}.
\]
According to \eqref{eq:matrix_model}, it holds that
\begin{align}
    \label{eq:exp_moment_matrix_model}
    &\notag
    \E \exp\left\{ \vec(\X)^\top V \vec(\X) - \Tr(V \Sigma) \right\}
    \\&\notag
    = \E \exp\left\{ \sum\limits_{j, k = 1}^K \vec(\Y_j)^\top (A_j^\top \otimes B_j^\top) V (A_k \otimes B_k) \vec(\Y_k) - \Tr(V \Sigma) \right\}
    \\&
    = \E \exp\Bigg\{ \sum\limits_{j = 1}^K \vec(\Y_j)^\top (A_j^\top \otimes B_j^\top) V (A_j \otimes B_j) \vec(\Y_j) - \Tr(V \Sigma)
    \\&\quad\notag
    + \sum\limits_{j \neq k} \vec(\Y_j)^\top (A_j^\top \otimes B_j^\top) V (A_k \otimes B_k) \vec(\Y_k) \Bigg\}.
\end{align}
Applying the Cauchy-Schwarz inequality, we bound the expression in the right-hand side of \eqref{eq:exp_moment_matrix_model} by a product of two terms:
\begin{align}
    \label{eq:exp_moment_cauchy-schwarz}
    &\notag
    \E \exp\left\{ \vec(\X)^\top V \vec(\X) - \Tr(V \Sigma) \right\}
    \\&
    \leq \left[ \E \exp\left\{ 2 \sum\limits_{j = 1}^K \vec(\Y_j)^\top (A_j^\top \otimes B_j^\top) V (A_j \otimes B_j) \vec(\Y_j) - 2 \Tr(V \Sigma) \right\} \right]^{1/2}
    \\& \quad \notag
    \cdot \left[ \E \exp\left\{ 2 \sum\limits_{j \neq k} \vec(\Y_j)^\top (A_j^\top \otimes B_j^\top) V (A_k \otimes B_k) \vec(\Y_k) \right\} \right]^{1/2}.
\end{align}
In the rest of the proof, we provide upper bounds on the exponential moments in the right-hand side of \eqref{eq:exp_moment_cauchy-schwarz}. For convenience, we split the proof into several steps.

\medskip

\noindent
\textbf{Step 1: an upper bound on the first term in \eqref{eq:exp_moment_cauchy-schwarz}.}
\quad
An upper bound on the first term in the right-hand side of \eqref{eq:exp_moment_cauchy-schwarz} simply follows from the conditions of the proposition. Since $\Y_1, \dots, \Y_K$ are independent, we have
\begin{align*}
    &
    \E \exp\left\{ 2 \sum\limits_{j = 1}^K \vec(\Y_j)^\top (A_j^\top \otimes B_j^\top) V (A_j \otimes B_j) \vec(\Y_j) - 2 \Tr(V \Sigma) \right\}
    \\&
    = \prod\limits_{j = 1}^K \E \exp\left\{ 2 \vec(\Y_j)^\top (A_j^\top \otimes B_j^\top) V (A_j \otimes B_j) \vec(\Y_j) - 2 \Tr\big((A_j^\top \otimes B_j^\top) V (A_j \otimes B_j) \big) \right\}.
\end{align*}
Then \eqref{eq:y_orlicz_norm} yields that
\begin{align}
    \label{eq:diagonal_exp_moment_bound}
    &\notag
    \E \exp\left\{ 2 \sum\limits_{j = 1}^K \vec(\Y_j)^\top (A_j^\top \otimes B_j^\top) V (A_j \otimes B_j) \vec(\Y_j) - 2 \Tr(V \Sigma) \right\}
    \\&
    \leq \prod\limits_{j = 1}^K \E \exp\left\{ 4 \omega_j^2 \left\|(A_j A_j^\top \otimes B_j B_j^\top)^{1/2} V (A_j A_j^\top \otimes B_j B_j^\top)^{1/2} \right\|_\Fr^2 \right\}
    \\&\notag
    \leq \exp\left\{ 4 \omega_{\max}^2 \sum\limits_{j = 1}^K \left\|(A_j A_j^\top \otimes B_j B_j^\top)^{1/2} V (A_j A_j^\top \otimes B_j B_j^\top)^{1/2} \right\|_\Fr^2 \right\}
\end{align}
for any $V$, such that
\[
    \left\|(A_j A_j^\top \otimes B_j B_j^\top)^{1/2} V (A_j A_j^\top \otimes B_j B_j^\top)^{1/2} \right\|_\Fr
    \leq \frac1{2 \omega_{\max}},
    \quad \text{where $\omega_{\max} = \max\limits_{1 \leq j \leq r} \omega_j$.}
\]

\medskip

\noindent
\textbf{Step 2: decoupling.}
\quad
We move to the analysis of the second term in the right-hand side of \eqref{eq:exp_moment_cauchy-schwarz}. The proof of an upper bound on
\[
    \E \exp\left\{ 2 \sum\limits_{j \neq k} \vec(\Y_j)^\top (A_j^\top \otimes B_j^\top) V (A_k \otimes B_k) \vec(\Y_k) \right\}
\]
extends the idea of decoupling (see \citep[Theorem 6.1.1]{vershynin18}) to the multivariate case. We proceed with the next auxiliary result.
\begin{Lem}
    \label{lem:decoupling}
    Let $\bzeta_1, \dots, \bzeta_K$ be i.i.d. centered random vectors in $\R^d$ and let $\{ M_{jk} : 1 \leq j \leq K, 1 \leq k \leq K\}$ be a collection of deterministic matrices of size $(d \times d)$. Then, for any convex function $G : \R \rightarrow \R$ with the finite expectation
    \[
        \E G \left( \sum\limits_{j \neq k} \bzeta_j^\top M_{jk} \bzeta_k \right),
    \]
    we have
    \[
        \E G \left( \sum\limits_{j \neq k} \bzeta_j^\top M_{jk} \bzeta_k \right) \leq \E G \left( 4 \sum\limits_{j \neq k} \bzeta_j^\top M_{jk} \bzeta'_k \right),
    \]
    where $\bzeta'_1, \dots, \bzeta'_K$ are independent copies of $\bzeta_1, \dots, \bzeta_K$. 
\end{Lem}
The proof of Lemma \ref{lem:decoupling} is postponed to Appendix \ref{sec:lem_decoupling_proof}. Applying this lemma to
\[
    \E \exp\left\{ 2 \sum\limits_{j \neq k} \vec(\Y_j)^\top (A_j^\top \otimes B_j^\top) V (A_k \otimes B_k) \vec(\Y_k) \right\},
\]
we obtain that
\begin{align}
    \label{eq:decoupling}
    &\notag
    \E \exp\left\{ 2 \sum\limits_{j \neq k} \vec(\Y_j)^\top (A_j^\top \otimes B_j^\top) V (A_k \otimes B_k) \vec(\Y_k) \right\}
    \\&
    \leq \E \exp\left\{ 8 \sum\limits_{j \neq k} \vec(\Y_j)^\top (A_j^\top \otimes B_j^\top) V (A_k \otimes B_k) \vec(\Y_k') \right\},
\end{align}
where $\Y_1', \dots, \Y_K'$ are independent copies of $\Y_1, \dots, \Y_K$. 

\medskip
\noindent\textbf{Step 3: reduction to the Gaussian case.}
\quad
Let $\bgamma_1, \dots, \bgamma_K \sim \N(\bzero, I_{d^2})$ be i.i.d Gaussian random vectors, which are independent of $\Y_1, \dots, \Y_K$.
We are going to show that the right-hand side of \eqref{eq:decoupling} does not exceed
\[
    \E \prod\limits_{j = 1}^K \exp\left\{ \frac{128 (1 + \omega_{\max})}{\log 2} \left\| \sum\limits_{k \neq j} (A_j^\top \otimes B_j^\top) V (A_k \otimes B_k) \bgamma_k\right\|^2 \right\}.
\]
In other words, it is enough to consider the Gaussian case. To do so, we first make a small detour and consider the $\psi_2$-norm of random variables of a form $\bu^\top \vec(\Y_j)$, where $\bu \in \R^{d^2}$ and $j \in \{1, \dots, K\}$.

\begin{Lem}
    \label{lem:psi_2_norm}
    Assume that $\Y_j$ satisfies \eqref{eq:y_orlicz_norm}, $j \in \{1, \dots, K\}$. Then, for any $\bu \in \R^{d^2}$, it holds that
    \[
        \left\| \bu^\top \vec(\Y_j) \right\|_{\psi_2}^2 \leq \frac{(1 + \omega_j) \|\bu\|^2}{\log 2}.
    \]
\end{Lem}
The proof of Lemma \ref{lem:psi_2_norm} is deferred to Appendix \ref{sec:lem_psi_2_norm_proof}. If a centered random variable has a finite $\psi_2$-norm, we can bound its exponential moments, using the following lemma.

\begin{Lem}
    \label{lem:orlicz_exp_moment}
    Let $\eta$ be a centered random variable with a finite Orlicz norm $\|\eta\|_{\psi_2} \leq \sigma < +\infty$. Then, for any $\lambda \in \R$, it holds that
    \[
        \E e^{\lambda \eta} \leq e^{\lambda^2 \sigma^2}.
    \]
\end{Lem}
It is well-known that random variables with a finite $\psi_2$-norm exhibit a sub-Gaussian behaviour (see, for instance, \citep[Section 2.6]{vershynin18}). However, the existing bounds in the literature include implicit absolute constants, while we are interested in tracking explicit ones. For this purpose, we carry out the proof of Lemma \ref{lem:orlicz_exp_moment} in Appendix \ref{sec:lem_orlicz_exp_moment_proof}. Using this lemma and considering the expectation
\[
    \E \left[ \exp\left\{ 8 \sum\limits_{j \neq k} \vec(\Y_j)^\top (A_j^\top \otimes B_j^\top) V (A_k \otimes B_k) \vec(\Y_k') \right\} \,\Bigg\vert\, \Y_1, \dots, \Y_K \right]
\]
conditionally on $\Y_1, \dots, \Y_K$, we obtain that
\begin{align}
    \label{eq:sub-gaussian_bound}
    &\notag
    \E \exp\left\{ 8 \sum\limits_{j \neq k} \vec(\Y_j)^\top (A_j^\top \otimes B_j^\top) V (A_k \otimes B_k) \vec(\Y_k') \right\}
    \\&
    \leq \prod\limits_{k = 1}^K \E \exp\left\{ \frac{64 (1 + \omega_{\max})}{\log 2} \left\| \sum\limits_{j \neq k}  (A_k^\top \otimes B_k^\top) V^\top (A_j \otimes B_j) \vec(\Y_j) \right\|^2 \right\}.
\end{align}
Now $\bgamma_1, \dots, \bgamma_K$ come into the play.
Using the fact that, for any $k \in \{1, \dots, K\}$ and any $\bu \in \R^{d^2}$, it holds that $\E e^{\bu^\top \bgamma_k} = e^{\|\bu\|^2 / 2}$, we represent the right-hand side of \eqref{eq:sub-gaussian_bound} in the following form:
\begin{align*}
    &
    \E_{\Y} \prod\limits_{k = 1}^K \exp\left\{ \frac{64 (1 + \omega_{\max})}{\log 2} \left\| \sum\limits_{j \neq k}  (A_k^\top \otimes B_k^\top) V^\top (A_j \otimes B_j) \vec(\Y_j) \right\|^2 \right\}
    \\&
    = \E_{\Y} \prod\limits_{k = 1}^K \E_{\bgamma_k} \exp\left\{ 8 \sqrt{\frac{2 + 2\omega_{\max}}{\log 2}} \cdot \bgamma_k^\top \sum\limits_{j \neq k}  (A_k^\top \otimes B_k^\top) V^\top (A_j \otimes B_j) \vec(\Y_j) \right\}
    \\&
    = \E \exp\left\{ 8 \sqrt{\frac{2 + 2\omega_{\max}}{\log 2}} \cdot \sum\limits_{j \neq k} \bgamma_k^\top (A_k^\top \otimes B_k^\top) V^\top (A_j \otimes B_j) \vec(\Y_j) \right\}.
\end{align*}
Applying Lemma \ref{lem:psi_2_norm} and Lemma \ref{lem:orlicz_exp_moment} again, we obtain that
\begin{align}
    \label{eq:gaussian_bound}
    &\notag
    \E \exp\left\{ 8 \sum\limits_{j \neq k} \vec(\Y_j)^\top (A_j^\top \otimes B_j^\top) V (A_k \otimes B_k) \vec(\Y_k') \right\}
    \\&
    \leq \E_\bgamma \prod\limits_{j = 1}^K \E_{\Y_j} \exp\left\{ 8 \sqrt{\frac{2 + 2\omega_{\max}}{\log 2}} \cdot \sum\limits_{k \neq j} \bgamma_k^\top (A_k^\top \otimes B_k^\top) V^\top (A_j \otimes B_j) \vec(\Y_j) \right\}
    \\&\notag
    \leq \E \prod\limits_{j = 1}^K \exp\left\{ \frac{128 (1 + \omega_{\max})}{\log 2} \left\| \sum\limits_{k \neq j} (A_j^\top \otimes B_j^\top) V (A_k \otimes B_k) \bgamma_k\right\|^2 \right\}.
\end{align}

\medskip

\noindent
\textbf{Step 4: exponential moments of a Gaussian quadratic form.}
\quad
Let $\bgamma = (\bgamma_1^\top, \dots, \bgamma_K^\top)^\top \in \R^{d^2 K}$ denote a standard Gaussian random vector, obtained by stacking $\bgamma_1, \dots, \bgamma_K$ together and let $S = S^\top \in \R^{d^2 K \times d^2 K}$ be such that
\begin{equation}
    \label{eq:quadratic_form}
    \frac{128 (1 + \omega_{\max})}{\log 2} \sum\limits_{j = 1}^K \left\| \sum\limits_{k \neq j} (A_j^\top \otimes B_j^\top) V (A_k \otimes B_k) \bgamma_k\right\|^2
    = \bgamma^\top S \bgamma.
\end{equation}
Obviously, $S \succeq O$, because the corresponding quadratic form is always non-negative. Lemma B.2 from \citep{spokoiny23} yields
\[
    \log \E e^{\bgamma^\top S \bgamma}
    \leq \Tr(S) + \frac{\Tr(S^2)}{1 - 2\|S\|}
    \leq \left(1 + \frac{\|S\|}{1 - 2\|S\|} \right) \Tr(S)
    \leq 2 \Tr(S),
\]
provided that $\|S\| \leq \Tr(S) \leq 1/4$. Note that
\begin{align}
    \label{eq:trace}
    &\notag
    \Tr(S)
    = \E \bgamma^\top S \bgamma
    \\&\notag
    = \frac{128 (1 + \omega_{\max})}{\log 2} \; \E \sum\limits_{j = 1}^K \left\| \sum\limits_{k \neq j} (A_j^\top \otimes B_j^\top) V (A_k \otimes B_k) \bgamma_k \right\|^2
    \\&
    = \frac{128 (1 + \omega_{\max})}{\log 2} \sum\limits_{j = 1}^K \left\| \sum\limits_{k \neq j} (A_j^\top \otimes B_j^\top) V (A_k \otimes B_k) \right\|_\Fr^2
    \\&\notag
    = \frac{128 (1 + \omega_{\max})}{\log 2} \left( \left\| \Sigma^{1/2} V \Sigma^{1/2} \right\|_{\Fr}^2 -  \sum\limits_{j = 1}^K \left\|(A_j A_j^\top \otimes B_j B_j^\top)^{1/2} V (A_j A_j^\top \otimes B_j B_j^\top)^{1/2} \right\|_\Fr^2 \right). 
\end{align}
Thus, it holds that
\begin{align}
    \label{eq:exp_moment_quadratic_form_bound}
    &
    \E \prod\limits_{j = 1}^K \exp\left\{ \frac{128 (1 + \omega_{\max})}{\log 2} \left\| \sum\limits_{k \neq j} (A_j^\top \otimes B_j^\top) V (A_k \otimes B_k) \bgamma_k\right\|^2 \right\}
    \\&\notag
    \leq \frac{256 (1 + \omega_{\max})}{\log 2} \left( \left\| \Sigma^{1/2} V \Sigma^{1/2} \right\|_{\Fr}^2 -  \sum\limits_{j = 1}^K \left\|(A_j A_j^\top \otimes B_j B_j^\top)^{1/2} V (A_j A_j^\top \otimes B_j B_j^\top)^{1/2} \right\|_\Fr^2 \right)
\end{align}
for all $V \in \R^{d^2 \times d^2}$, such that
\[
    \frac{128 (1 + \omega_{\max})}{\log 2} \left\| \Sigma^{1/2} V \Sigma^{1/2} \right\|_{\Fr}^2 \leq \frac14.
\]
Hence, due to the inequalities \eqref{eq:decoupling}, \eqref{eq:gaussian_bound}, and \eqref{eq:exp_moment_quadratic_form_bound}, we have
\begin{align}
    \label{eq:off-diagonal_exp_moment_bound}
    &
    \E \exp\left\{ 2 \sum\limits_{j \neq k} \vec(\Y_j)^\top (A_j^\top \otimes B_j^\top) V (A_k \otimes B_k) \vec(\Y_k) \right\}
    \\&\notag \hspace{-7pt}
    \leq \exp\left\{256 (1 + \omega_{\max}) \left( \left\| \Sigma^{1/2} V \Sigma^{1/2} \right\|_{\Fr}^2 -  \sum\limits_{j = 1}^K \left\|(A_j A_j^\top \otimes B_j B_j^\top)^{1/2} V (A_j A_j^\top \otimes B_j B_j^\top)^{1/2} \right\|_\Fr^2 \right) \right\},
\end{align}
provided that
\[
    \frac{128 (1 + \omega_{\max})}{\log 2} \left\| \Sigma^{1/2} V \Sigma^{1/2} \right\|_{\Fr}^2 \leq \frac14.
\]

\medskip

\noindent
\textbf{Step 5: final bound.}
\quad
Taking into account \eqref{eq:diagonal_exp_moment_bound} and \eqref{eq:off-diagonal_exp_moment_bound}, we obtain that
\begin{align*}
    &
    \E \exp\left\{ \vec(\X)^\top V \vec(\X) - \Tr(V \Sigma) \right\}
    \\&
    \leq \left[ \E \exp\left\{ 2 \sum\limits_{j = 1}^K \vec(\Y_j)^\top (A_j^\top \otimes B_j^\top) V (A_j \otimes B_j) \vec(\Y_j) - 2 \Tr(V \Sigma) \right\} \right]^{1/2}
    \\& \:
    \cdot \left[ \E \exp\left\{ 2 \sum\limits_{j \neq k} \vec(\Y_j)^\top (A_j^\top \otimes B_j^\top) V (A_k \otimes B_k) \vec(\Y_k) \right\} \right]^{1/2}
    \\&
    \leq \exp\left\{ 2 \omega_{\max}^2 \sum\limits_{j = 1}^K \left\|(A_j A_j^\top \otimes B_j B_j^\top)^{1/2} V (A_j A_j^\top \otimes B_j B_j^\top)^{1/2} \right\|_\Fr^2 \right\}
    \\& \:
    \cdot \exp\left\{ \frac{128 (1 + \omega_{\max})}{\log 2} \left( \left\| \Sigma^{1/2} V \Sigma^{1/2} \right\|_{\Fr}^2 - \sum\limits_{j = 1}^K \left\|(A_j A_j^\top \otimes B_j B_j^\top)^{1/2} V (A_j A_j^\top \otimes B_j B_j^\top)^{1/2} \right\|_\Fr^2 \right) \right\}
    \\&
    \leq \exp\left\{ \left( \frac{128 (1 + \omega_{\max})}{\log 2} \vee 2 \omega_{\max}^2 \right) \left\| \Sigma^{1/2} V \Sigma^{1/2} \right\|_{\Fr}^2 \right\},
\end{align*}
provided that
\begin{equation}
    \label{eq:j-th_term_condition}
    \left\|(A_j A_j^\top \otimes B_j B_j^\top)^{1/2} V (A_j A_j^\top \otimes B_j B_j^\top)^{1/2} \right\|_\Fr
    \leq \frac1{2 \omega_{\max}}
\end{equation}
and
\[
    \frac{128 (1 + \omega_{\max})}{\log 2} \left\| \Sigma^{1/2} V \Sigma^{1/2} \right\|_{\Fr}^2 \leq \frac14.
\]
Note that, according to \eqref{eq:trace},
\[
    \left\| \Sigma^{1/2} V \Sigma^{1/2} \right\|_{\Fr}^2 -  \sum\limits_{j = 1}^K \left\|(A_j A_j^\top \otimes B_j B_j^\top)^{1/2} V (A_j A_j^\top \otimes B_j B_j^\top)^{1/2} \right\|_\Fr^2
    = \frac{\Tr(S) \, \log 2}{128 (1 + \omega_{\max})} 
    \geq 0,
\]
where $S$ is a positive semidefinite matrix, defined in \eqref{eq:quadratic_form}. Thus, the condition
\[
    \left\| \Sigma^{1/2} V \Sigma^{1/2} \right\|_{\Fr}
    \leq \frac1{2 \omega_{\max}}
\]
automatically yields \eqref{eq:j-th_term_condition}. Hence, we proved that
\[
    \E \exp\left\{ \vec(\X)^\top V \vec(\X) - \Tr(V \Sigma) \right\}
    \leq
    \exp\left\{ \left( \frac{128 (1 + \omega_{\max})}{\log 2} \vee 2 \omega_{\max}^2 \right) \left\| \Sigma^{1/2} V \Sigma^{1/2} \right\|_{\Fr}^2 \right\}
\]
for all $V \in \R^{d^2 \times d^2}$, such that
\[
    \left(4 \omega_{\max}^2 \vee \frac{512 (1 + \omega_{\max})}{\log 2} \right) \left\| \Sigma^{1/2} V \Sigma^{1/2} \right\|_{\Fr}^2 \leq 1.
\]
\endproof

\section{Proof of Lemma \ref{lem:ro_representation}}
\label{sec:lem_ro_representation_proof}

We start with representing the operator norm of $\ro(\widehat \Sigma - \Sigma)$ in the following form:
\[
    \left\| \ro(\widehat \Sigma - \Sigma) \right\|
    = \sup\limits_{\substack{U \in \R^{p \times p}, V \in \R^{q \times q}\\ \|U\|_\Fr = \|V\|_\Fr = 1}} \vec(U)^\top \ro(\widehat \Sigma - \Sigma) \vec(V).
\]
Due to the linearity of the rearrangement operator $\ro$, we have
\begin{align}
    \label{eq:ro_supremum}
    \left\| \ro(\widehat \Sigma - \Sigma) \right\|
    &
    = \sup\limits_{\substack{U \in \R^{p \times p}, V \in \R^{q \times q}\\ \|U\|_\Fr = \|V\|_\Fr = 1}} \vec(U)^\top \ro\left( \frac1n \sum\limits_{i = 1}^n \bX_i \bX_i^\top - \E \bX \bX^\top \right) \vec(V)
    \\&\notag
    = \sup\limits_{\substack{U \in \R^{p \times p}, V \in \R^{q \times q}\\ \|U\|_\Fr = \|V\|_\Fr = 1}} \left\{ \frac1n \sum\limits_{i = 1}^n \vec(U)^\top \left[ \ro\left( \bX_i \bX_i^\top \right) - \E \ro\left(\bX \bX^\top \right) \right] \vec(V) \right\}.
\end{align}
Let us show that, for any $U \in \R^{p \times p}$ and $V \in \R^{q \times q}$,
\[
    \vec(U)^\top \ro\left(\bX \bX^\top \right) \vec(V)
    = \Tr\left( \X^\top V^\top \X U \right),
\]
where $\X$ is a matrix of size $(q \times p)$, such that $\vec(\X) = \bX$. It holds that
\begin{align*}
    &
    \vec(U)^\top \ro(\bX \bX^\top) \vec(V)
    \\&
    = \sum\limits_{i = 1}^p \sum\limits_{j = 1}^p \sum\limits_{k = 1}^q \sum\limits_{\ell = 1}^q \vec(U)_{p (i - 1) + j} \ro(\bX \bX^\top)_{p (i - 1) + j, q (k - 1) + \ell} \vec(V)_{q (k - 1) + \ell}.
\end{align*}
The definition of the vectorization operator $\vec$ yields that
\[
    \vec(U)_{p (i - 1) + j} = U_{ji}
    \quad \text{and} \quad
    \vec(V)_{q (k - 1) + \ell} = V_{\ell k}.
\]
At the same time, according to the definition of the rearrangement operator $\ro$, we have
\begin{align*}
    \ro(\bX \bX^\top)_{p (i - 1) + j, q (k - 1) + \ell}
    &
    = \vec\left( (\bX \bX^\top)_{p (j - 1) : pj, p (i - 1) : pi} \right)_{q (k - 1) + \ell}
    \\&
    = \left( (\bX \bX)^\top_{(p - 1) j : pj, (p - 1) i : pi} \right)_{\ell k}
    \\&
    = (\bX \bX)^\top_{(p - 1) j + \ell, (p - 1) i + k}
    \\&
    = \bX_{(p - 1) j + \ell} \bX_{(p - 1) i + k}
    = \X_{\ell j} \X_{k i}.
\end{align*}
Here $(\bX \bX^\top)_{p (j - 1) : pj, p (i - 1) : pi}$ denotes a submatrix of $\bX \bX^\top$ on the intersection of the rows $p (j - 1) + 1, \dots, pj$ with the columns $p (i - 1) + 1, \dots, pi$.
Hence, we obtain that
\begin{equation}
    \label{eq:ro_x_reshaping}
    \vec(U)^\top \ro(\bX \bX^\top) \vec(V)
    = \sum\limits_{i = 1}^p \sum\limits_{j = 1}^p \sum\limits_{k = 1}^q \sum\limits_{\ell = 1}^q U_{ji} \X_{\ell j} \X_{k i} V_{\ell k}
    = \Tr(V^\top \X U \X^\top).
\end{equation}
Similarly, one can prove that
\begin{equation}
    \label{eq:ro_x_i_reshaping}
    \vec(U)^\top \ro\left(\bX_i \bX_i^\top \right) \vec(V)
    = \Tr\left( \X_i^\top V^\top \X_i U \right),
\end{equation}
for any $U \in \R^{p \times p}$, $V \in \R^{q \times q}$, and any $i \in \{1, \dots, n\}$, where $\X_1, \dots, \X_n$ are the matrices of size $(q \times p)$, obtained by reshaping $\bX_1, \dots, \bX_n$, respectively:
\[
    \vec(\X_i) = \bX_i
    \quad \text{for all $i \in \{1, \dots, n\}$.}
\]
Taking into account \eqref{eq:ro_supremum}, \eqref{eq:ro_x_reshaping}, and \eqref{eq:ro_x_i_reshaping}, we deduce that
\begin{align*}
    \left\| \ro(\widehat \Sigma - \Sigma) \right\|
    &
    = \sup\limits_{\substack{U \in \R^{p \times p}, V \in \R^{q \times q}\\ \|U\|_\Fr = \|V\|_\Fr = 1}} \left\{ \frac1n \sum\limits_{i = 1}^n \vec(U)^\top \left[ \ro\left( \bX_i \bX_i^\top \right) - \E \ro\left(\bX \bX^\top \right) \right] \vec(V) \right\}
    \\&
    = \sup\limits_{\substack{U \in \R^{p \times p}, V \in \R^{q \times q}\\ \|U\|_\Fr = \|V\|_\Fr = 1}} \left\{ \frac1n \sum\limits_{i = 1}^n \Tr\left( \X_i^\top V^\top \X_i U \right) - \E \Tr\left( \X^\top V^\top \X U \right) \right\}.
\end{align*}

\endproof

\section{Proof of Lemma \ref{lem:z}}
\label{sec:lem_z_proof}
Let $\X \in \R^{p \times p}$ and $\Y \in \R^{q \times q}$ be random matrices with the joint distribution $\mu$, defined in \eqref{eq:prior_density}. That is, the entries of $\X$ and $\Y$ are independent centered Gaussian random variables, such that
\[
    \Var(\X_{ij}) = \frac1\alpha,
    \quad
    \Var(\Y_{k \ell}) = \frac1\beta
    \quad \text{for all $i, j \in \{1, \dots, p\}$ and $k, \ell \in \{1, \dots, q\}$.}
\]
We are going to show that
\begin{equation}
    \label{eq:fr_norm_1_bound}
    \p\left( \left\| \Sigma^{1/2} (\X \otimes \Y) \Sigma^{1/2} \right\|_\Fr^2 \geq \frac{4 \Tr(\Sigma)^2}{\alpha \beta} \right)
    \leq \frac14, 
\end{equation}
\begin{equation}
    \label{eq:fr_norm_2_bound}
    \p\left( \left\| \left( \sum\limits_{j = 1}^K \|\Psi_j\| \Phi_j \right)^{1/2} \X \left( \sum\limits_{j = 1}^K \|\Psi_j\| \Phi_j \right)^{1/2} \right\|_{\Fr}^2 \geq \frac{4}{\alpha} \left( \sum\limits_{j = 1}^K \|\Psi_j\| \Tr(\Phi_j) \right)^{2} \right)
    \leq \frac14,
\end{equation}
and
\begin{equation}
    \label{eq:fr_norm_3_bound}
    \p\left( \left\| \left( \sum\limits_{j = 1}^K \|\Phi_j\| \Psi_j \right)^{1/2} \Y \left( \sum\limits_{j = 1}^K \|\Phi_j\| \Psi_j \right)^{1/2} \right\|_{\Fr}^2 \leq \frac{4}{\beta} \left( \sum\limits_{j = 1}^K \|\Phi_j\| \Tr(\Psi_j) \right)^{2} \right)
    \leq \frac14.
\end{equation}
Then the union bound implies that
\[
    \ttZ
    = \p\left( (\X, \Y) \in \Upsilon(O, O) \right)
    = 1 - \p\left( (\X, \Y) \notin \Upsilon(O, O) \right)
    \geq 1 - \frac34
    = \frac14.
\]
For convenience, the rest of the proof is divided into three parts.

\medskip

\noindent
\textbf{Step 1: proof of \eqref{eq:fr_norm_1_bound}.}\quad
Since the covariance matrix $\Sigma$ has a Kronecker product structure
\[
    \Sigma = \sum\limits_{j = 1}^K \Phi_j \otimes \Psi_j,
\]
it holds that
\begin{align*}
    \left\| \Sigma^{1/2} (\X \otimes \Y) \Sigma^{1/2} \right\|_\Fr^2
    &
    = \Tr\big( \Sigma (\X \otimes \Y)^\top \Sigma (\X \otimes \Y) \big)
    \\&
    = \sum\limits_{j = 1}^K \sum\limits_{k = 1}^K \Tr\left( (\Phi_j \otimes \Psi_j) (\X^\top \otimes \Y^\top) (\Phi_k \otimes \Psi_k) (\X \otimes \Y) \right)
    \\&
    = \sum\limits_{j = 1}^K \sum\limits_{k = 1}^K \Tr\left( (\Phi_j \X^\top \Phi_k \X) \otimes (\Psi_j \Y^\top \Psi_k \Y) \right).
\end{align*}
Using the property \eqref{eq:kronecker_eigenvalues} about the trace of Kronecker product, we obtain that
\begin{align*}
    \left\| \Sigma^{1/2} (\X \otimes \Y) \Sigma^{1/2} \right\|_\Fr^2
    &
    = \sum\limits_{j = 1}^K \sum\limits_{k = 1}^K \Tr\left( (\Phi_j \X^\top \Phi_k \X) \otimes (\Psi_j \Y^\top \Psi_k \Y) \right)
    \\&
    = \sum\limits_{j = 1}^K \sum\limits_{k = 1}^K \Tr(\Phi_j \X^\top \Phi_k \X) \; \Tr(\Psi_j \Y^\top \Psi_k \Y).
\end{align*}
For any $j$ and $k$ from $\{1, \dots, K\}$, we have
\[
    \E \Tr(\Phi_j \X^\top \Phi_k \X)
    = \E \vec(\X)^\top (\Phi_j \otimes \Phi_k) \vec(\X)
    = \frac1\alpha \Tr(\Phi_j \otimes \Phi_k)
    = \frac1\alpha \; \Tr(\Phi_j) \; \Tr(\Phi_k)
\]
and, similarly,
\[
    \E \Tr(\Psi_j \Y^\top \Psi_k \Y)
    = \E \vec(\Y)^\top (\Psi_j \otimes \Psi_k) \vec(\Y)
    = \frac1\beta \Tr(\Psi_j \otimes \Psi_k)
    = \frac1\beta \; \Tr(\Psi_j) \; \Tr(\Psi_k).
\]
This yields that
\begin{align*}
    \E \left\| \Sigma^{1/2} (\X \otimes \Y) \Sigma^{1/2} \right\|_\Fr^2
    &
    = \frac1{\alpha \beta} \sum\limits_{j = 1}^K \sum\limits_{k = 1}^K \Tr(\Phi_j) \Tr(\Phi_k) \Tr(\Psi_j) \Tr(\Psi_k)
    \\&
    = \frac1{\alpha \beta} \left( \sum\limits_{j = 1}^K \Tr(\Phi_j) \Tr(\Psi_j)\right)^2
    = \frac{\Tr(\Sigma)^2}{\alpha \beta},
\end{align*}
and, applying Markov's inequality, we obtain \eqref{eq:fr_norm_1_bound}.

\medskip

\noindent
\textbf{Step 2: proof of \eqref{eq:fr_norm_2_bound} and \eqref{eq:fr_norm_3_bound}.}\quad
Similarly to the first step, the proof of \eqref{eq:fr_norm_2_bound} and \eqref{eq:fr_norm_3_bound} relies on Markov's inequality and computation of the expectations
\[
    \E \left\| \left( \sum\limits_{j = 1}^K \|\Psi_j\| \Phi_j \right)^{1/2} \X \left( \sum\limits_{j = 1}^K \|\Psi_j\| \Phi_j \right)^{1/2} \right\|_{\Fr}^2
\]
and
\[
    \E \left\| \left( \sum\limits_{j = 1}^K \|\Phi_j\| \Psi_j \right)^{1/2} \Y \left( \sum\limits_{j = 1}^K \|\Phi_j\| \Psi_j \right)^{1/2} \right\|_{\Fr}^2.
\]
First, let us show that
\[
    \E \left\| \left( \sum\limits_{j = 1}^K \|\Psi_j\| \Phi_j \right)^{1/2} \X \left( \sum\limits_{j = 1}^K \|\Psi_j\| \Phi_j \right)^{1/2} \right\|_{\Fr}^2
    = \frac{1}{\alpha} \left( \sum\limits_{j = 1}^K \|\Psi_j\| \; \Tr(\Phi_j) \right)^{2}.
\]
Indeed, it holds that
\begin{align*}
    &
    \E \left\| \left( \sum\limits_{j = 1}^K \|\Psi_j\| \Phi_j \right)^{1/2} \X \left( \sum\limits_{j = 1}^K \|\Psi_j\| \Phi_j \right)^{1/2} \right\|_{\Fr}^2
    \\&
    = \E \Tr \left[ \left( \sum\limits_{j = 1}^K \|\Psi_j\| \Phi_j \right) \X^\top \left( \sum\limits_{k = 1}^K \|\Psi_k\| \Phi_k \right) \X \right]
    \\&
    = \sum\limits_{j = 1}^K \sum\limits_{k = 1}^K \|\Psi_j\| \|\Psi_k\| \; \E \Tr \left(\Phi_j \X^\top \Phi_k \X \right).
\end{align*}
Since for any $j$ and $k$ from $\{1, \dots, K\}$ the expectation of $\Tr \left(\Phi_j \X^\top \Phi_k \X \right)$ is equal to
\[
    \E \Tr \left(\Phi_j \X^\top \Phi_k \X \right)
    = \E \vec(\X)^\top (\Phi_j \otimes \Phi_k) \vec(\X)
    = \frac1{\alpha} \Tr(\Phi_j \otimes \Phi_k)
    = \frac1{\alpha} \Tr(\Phi_j) \; \Tr(\Phi_k),
\]
we obtain that
\begin{align*}
    \E \left\| \left( \sum\limits_{j = 1}^K \|\Psi_j\| \Phi_j \right)^{1/2} \X \left( \sum\limits_{j = 1}^K \|\Psi_j\| \Phi_j \right)^{1/2} \right\|_{\Fr}^2
    &
    = \frac1{\alpha} \sum\limits_{j = 1}^K \sum\limits_{k = 1}^K \|\Psi_j\| \|\Psi_k\| \; \Tr(\Phi_j) \Tr(\Phi_k)
    \\&
    = \frac1{\alpha} \left( \sum\limits_{j = 1}^K \|\Psi_j\| \; \Tr(\Phi_j) \right)^2,
\end{align*}
and \eqref{eq:fr_norm_2_bound} follows from Markov's inequality. The proof of \eqref{eq:fr_norm_3_bound} is absolutely similar.

\endproof

\section{Proof of Lemma \ref{lem:support}}
\label{sec:lem_support_proof}

Note that
\begin{align*}
    \left\|\Sigma^{1/2} (P \otimes Q) \Sigma^{1/2} \right\|_{\Fr}^2
    &
    \leq 4 \left\|\Sigma^{1/2} \big( (P - U) \otimes (Q - V) \big) \Sigma^{1/2} \right\|_{\Fr}^2
    \\&\quad
    + 4 \left\|\Sigma^{1/2} \big(U \otimes (Q - V) \big) \Sigma^{1/2} \right\|_{\Fr}^2
    \\&\quad
    + 4 \left\|\Sigma^{1/2} \big((P - U) \otimes V \big) \Sigma^{1/2} \right\|_{\Fr}^2
    \\&\quad
    + 4 \left\|\Sigma^{1/2} (U \otimes V) \Sigma^{1/2} \right\|_{\Fr}^2.
\end{align*}
The definition of the support $\Upsilon(U, V)$ (see \eqref{eq:support}) implies that
\begin{equation}
    \label{eq:support_1_condition}
    \|\Sigma^{1/2} \big( (P - U) \otimes (Q - V) \big) \Sigma^{1/2}\|_{\Fr}^2
    \leq \frac{4 \Tr(\Sigma)^2}{\alpha \beta}
    \quad
    \text{$\rho_{U,V}$-almost surely.}
\end{equation}
Moreover, on $\Upsilon(U, V)$ it holds that
\begin{align*}
    &
    \left\|\Sigma^{1/2} \big((P - U) \otimes V \big) \Sigma^{1/2} \right\|_{\Fr}^2
    = \Tr\left(\Sigma \big((P - U) \otimes V \big)^\top \Sigma \big((P - U) \otimes V \big) \right)
    \\&
    = \sum\limits_{j = 1}^K \sum\limits_{k = 1}^K \Tr\left((\Phi_j \otimes \Psi_j) \big((P - U) \otimes V \big)^\top (\Phi_k \otimes \Psi_k) \big((P - U) \otimes V \big) \right)
    \\&
    = \sum\limits_{j = 1}^K \sum\limits_{k = 1}^K \Tr\big(\Phi_j (P - U)^\top \Phi_k (P - U) \big) \Tr\big(\Psi_j V^\top \Psi_k V) \big)
    \\&
    = \sum\limits_{j = 1}^K \sum\limits_{k = 1}^K \Tr\big(\Phi_j (P - U)^\top \Phi_k (P - U) \big) \left\| \Psi_j^{1/2} V \Psi_k^{1/2} \right\|_{\Fr}^2.
\end{align*}
Since $\|V\|_\Fr = 1$, we have
\[
    \left\| \Psi_j^{1/2} V \Psi_k^{1/2} \right\|_{\Fr}^2
    \leq \left\| \Psi_j^{1/2} \right\|^2 \cdot \big\|V \big\|_\Fr^2 \cdot \left\| \Psi_k^{1/2} \right\|^2
    \leq \|\Psi_j\| \|\Psi_k\|
    \quad \text{for all $j, k \in \{1, \dots, K\}$,}
\]
and thus,
\begin{align}
    \label{eq:support_2_condition}
    \left\|\Sigma^{1/2} \big((P - U) \otimes V \big) \Sigma^{1/2} \right\|_{\Fr}^2
    &\notag
    \leq \sum\limits_{j = 1}^K \sum\limits_{k = 1}^K \|\Psi_j\| \|\Psi_k\| \Tr\big(\Phi_j (P - U)^\top \Phi_k (P - U) \big)
    \\&
    = \left\| \left(\sum\limits_{j = 1}^K \|\Psi_j\| \Phi_j\right)^{1/2} (P - U) \left(\sum\limits_{k = 1}^K \|\Psi_k\| \Phi_k\right)^{1/2} \right\|_\Fr^2
    \\&\notag
    \leq \frac{4}{\alpha} \left( \sum\limits_{j = 1}^K \|\Psi_j\| \Tr(\Phi_j) \right)^{2}
    \quad
    \text{$\rho_{U,V}$-almost surely,}
\end{align}
where the last inequality is due to \eqref{eq:support}. Similarly, one can show that
\begin{equation}
    \label{eq:support_3_condition}
    \left\|\Sigma^{1/2} \big(U \otimes (Q - V) \big) \Sigma^{1/2} \right\|_{\Fr}^2
    \leq \frac{4}{\beta} \left( \sum\limits_{j = 1}^K \|\Phi_j\| \Tr(\Psi_j) \right)^{2}
    \quad
    \text{$\rho_{U,V}$-almost surely.}
\end{equation}
Finally,
\[
    \left\|\Sigma^{1/2} (U \otimes V) \Sigma^{1/2} \right\|_{\Fr}^2
    \leq \left\|\Sigma^{1/2} \right\|^2 \cdot \left\|U \otimes V \right\|_\Fr^2 \cdot \left\|\Sigma^{1/2} \right\|^2
    = \left\|\Sigma \right\|^2.
\]
This inequality and the bounds \eqref{eq:support_1_condition}, \eqref{eq:support_2_condition}, \eqref{eq:support_3_condition} yield that
\begin{align*}
    \left\|\Sigma^{1/2} (P \otimes Q) \Sigma^{1/2} \right\|_{\Fr}^2
    &
    \leq \left\|\Sigma \right\|^2 + \frac{4}{\alpha} \left( \sum\limits_{j = 1}^K \|\Psi_j\| \Tr(\Phi_j) \right)^{2}
    \\&\quad
    + \frac{4}{\beta} \left( \sum\limits_{j = 1}^K \|\Phi_j\| \Tr(\Psi_j) \right)^{2} + \frac{4 \Tr(\Sigma)^2}{\alpha \beta}
    \quad \text{$\rho_{U,V}$-almost surely.}
\end{align*}
It only remains to recall the definition of $\alpha$ and $\beta$ (see \eqref{eq:alpha_beta}) to conclude the proof:
\begin{align*}
    &
    \left\|\Sigma^{1/2} (P \otimes Q) \Sigma^{1/2} \right\|_{\Fr}^2
    \\&
    \leq \left\|\sum\limits_{j = 1}^K \Phi_j \otimes \Psi_j \right\|^2 + \frac{4}{\max\limits_{1 \leq j \leq K} \ttr(\Phi_j)^2} \left( \sum\limits_{j = 1}^K \|\Psi_j\| \Tr(\Phi_j) \right)^{2}
    \\&\quad
    + \frac{4}{\max\limits_{1 \leq j \leq K} \ttr(\Psi_j)^2}  \left( \sum\limits_{j = 1}^K \|\Phi_j\| \Tr(\Psi_j) \right)^{2} + \frac{4}{\max\limits_{1 \leq j \leq K} \ttr(\Phi_j)^2 \ttr(\Psi_j)^2} \left( \sum\limits_{j = 1}^K \Tr(\Phi_j) \Tr(\Psi_j) \right)^2
    \\&
    \leq \sum\limits_{j = 1}^K \left\| \Phi_j \otimes \Psi_j \right\|^2 + 12 \left( \sum\limits_{j = 1}^K \|\Phi_j\| \|\Psi_j\| \right)^2
    \\&
    = 13 \left( \sum\limits_{j = 1}^K \|\Phi_j\| \|\Psi_j\| \right)^2
    \quad \text{$\rho_{U,V}$-almost surely.}
\end{align*}
\endproof

\section{Proof of Lemma \ref{lem:decoupling}}
\label{sec:lem_decoupling_proof}

Let $\eps_1, \dots, \eps_K \sim \Be(1/2)$ be i.i.d. Bernoulli random variables. Then
\[
    \sum\limits_{j \neq k} \bzeta_j^\top M_{jk} \bzeta_k
    = 4 \E_\eps \sum\limits_{j \neq k} \eps_j (1 - \eps_k) \bzeta_j^\top M_{jk} \bzeta_k
    \quad \text{almost surely,}
\]
where $\E_\eps$ stands for the expectation with respect to $\eps_1, \dots, \eps_K$ conditionally on $\bzeta_1, \dots, \bzeta_K$. Let us define a random set $\I \subseteq \{1, \dots, K\}$ as follows:
\[
    \I = \left\{ j \in \{1, \dots, K\} : \eps_j = 1 \right\}.
\]
Then we can rewrite the sum
\[
    \sum\limits_{j \neq k} \eps_j (1 - \eps_k) \bzeta_j^\top M_{jk} \bzeta_k
\]
in the form
\[
    \sum\limits_{j \neq k} \eps_j (1 - \eps_k) \bzeta_j^\top M_{jk} \bzeta_k
    = \sum\limits_{j \in \I} \sum\limits_{k \in \I^c} \bzeta_j^\top M_{jk} \bzeta_k,
\]
where $\I^c = \{1, \dots, K\} \backslash \I$ denotes the complement of $\I$. With the notation introduced above, we obtain that
\[
    \E G \left( \sum\limits_{j \neq k} \bzeta_j^\top M_{jk} \bzeta_k \right)
    = \E_{\bzeta} G \left( 4 \E_\eps \sum\limits_{j \in \I} \sum\limits_{k \in \I^c} \bzeta_j^\top M_{jk} \bzeta_k \right).
\]
Here $\E_{\bzeta}$ stands for the expectation with respect to $\bzeta_1, \dots, \bzeta_K$. Applying Jensen's inequality and taking into account that, for any realization of $\eps_1, \dots, \eps_K$, the collection $\{\bzeta_j : j \in \I\}$ is independent of $\{\bzeta_k : k \in \I\}$, we deduce that
\[
    \E G \left( \sum\limits_{j \neq k} \bzeta_j^\top M_{jk} \bzeta_k \right)
    \leq \E_\eps \E_{\bzeta} G \left( 4 \sum\limits_{j \in \I} \sum\limits_{k \in \I^c} \bzeta_j^\top M_{jk} \bzeta_k \right)
    \leq \E_\eps \E_{\bzeta, \bzeta'} G \left( 4 \sum\limits_{j \in \I} \sum\limits_{k \in \I^c} \bzeta_j^\top M_{jk} \bzeta_k' \right).
\]
Note that, for any realization of $\eps_1, \dots, \eps_r$, the sum
\[
    \sum\limits_{j \in \I^c} \sum\limits_{k \in \I} \bzeta_j^\top M_{jk} \bzeta_k'
    \quad \text{is independent of} \quad
    \sum\limits_{j \in \I} \sum\limits_{k \in \I^c} \bzeta_j^\top M_{jk} \bzeta_k'
\]
and
\[
   \E \sum\limits_{j \in \I^c} \sum\limits_{k \in \I} \bzeta_j^\top M_{jk} \bzeta_k' = 0. 
\]
Then Lemma 6.1.2 from the book \citep{vershynin18} implies that
\begin{align*}
    \E_{\bzeta, \bzeta'} G \left( 4 \sum\limits_{j \in \I} \sum\limits_{k \in \I^c} \bzeta_j^\top M_{jk} \bzeta_k' \right)
    &
    \leq \E_{\bzeta, \bzeta'} G \left( 4 \sum\limits_{j \in \I} \sum\limits_{k \in \I^c} \bzeta_j^\top M_{jk} \bzeta_k' + 4 \sum\limits_{j \in \I^c} \sum\limits_{k \in \I} \bzeta_j^\top M_{jk} \bzeta_k' \right)
    \\&
    = \E_{\bzeta, \bzeta'} G \left( 4 \sum\limits_{j \neq k} \bzeta_j^\top M_{jk} \bzeta_k' \right).
\end{align*}
Hence, we obtain that
\begin{align*}
    \E G \left( \sum\limits_{j \neq k} \bzeta_j^\top M_{jk} \bzeta_k \right)
    &
    \leq \E_\eps \E_{\bzeta, \bzeta'} G \left( 4 \sum\limits_{j \in \I} \sum\limits_{k \in \I^c} \bzeta_j^\top M_{jk} \bzeta_k' \right)
    \\&
    \leq \E_\eps \E_{\bzeta, \bzeta'} G \left( 4 \sum\limits_{j \neq k} \bzeta_j^\top M_{jk} \bzeta_k' \right)
    \\&
    = \E_{\bzeta, \bzeta'} G \left( 4 \sum\limits_{j \neq k} \bzeta_j^\top M_{jk} \bzeta_k' \right).
\end{align*}
\endproof

\section{Proof of Lemma \ref{lem:psi_2_norm}}
\label{sec:lem_psi_2_norm_proof}
Let us denote
\[
    t = \|\bu\| \sqrt{\frac{1 + \omega_j}{\log 2}}
    \geq \|\bu\| \sqrt{\omega_j} 
\]
and show that
\[
    \E \exp\left\{ \frac{\big(\bu^\top \vec(\Y_j) \big)^2}{t^2} \right\} \leq 2.
\]
Since $\vec(\Y_j)$ satisfies \eqref{eq:y_orlicz_norm}, we have
\begin{align}
    \label{eq:psi_2_norm_bound}
    \log \E \exp\left\{ \frac{\big(\bu^\top \vec(\Y_j) \big)^2}{t^2} \right\}
    &\notag
    \leq \frac{\Tr(\bu \bu^\top)}{t^2} + \frac{\omega_j^2 \left\|\bu \bu^\top \right\|_{\Fr}^2}{t^4}
    \\&
    = \frac{\|\bu\|^2}{t^2} + \frac{\omega_j^2 \|\bu\|^4}{t^4}
    = \frac{\|\bu\|^2}{t^2} \left( 1 + \omega_j \cdot \frac{\omega_j \|\bu\|^2}{t^2} \right)
    \\&\notag
    \leq \frac{(1 + \omega_j) \|\bu\|^2}{t^2}
    \leq \log 2.
\end{align}
Here we took into account that $\omega_j \|\bu\|^2 < t^2$. The inequality \eqref{eq:psi_2_norm_bound} yields that
\[
    \left\| \bu^\top \vec(\Y_j) \right\|_{\psi_2}^2 \leq \frac{(1 + \omega_j) \|\bu\|^2}{\log 2}.
\]
\endproof

\section{Proof of Lemma \ref{lem:orlicz_exp_moment}}
\label{sec:lem_orlicz_exp_moment_proof}

\noindent
\textbf{Step 1: a bound on $\E \eta^{2k}$.}
\quad
The goal of this step is to prove that
\[
    \E \eta^{2k} \leq 2 \cdot k! \cdot \sigma^{2k}.
\]
Note that Markov's inequality and the definition of the Orlicz norm yield that
\[
    \p\left( |\eta| \geq t \right)
    \leq e^{-t^2 / \sigma^2} \E e^{\eta^2 / \sigma^2}
    \leq 2 e^{-t^2 / \sigma^2}
    \quad \text{for all $t > 0$.}
\]
Then it holds that
\[
    \E \eta^{2k}
    = \int\limits_0^{+\infty} \p\left( \eta^{2k} \geq t \right) \dd t
    = \int\limits_0^{+\infty} \p\left( |\eta| \geq t^{1 / (2k)} \right) \dd t
    \leq 2 \int\limits_0^{+\infty} \exp\left\{ -\frac{t^{1/k}}{\sigma^2} \right\} \dd t.
\]
Substituting $t$ with $\sigma^{2k} u^k$, we obtain that
\begin{equation}
    \label{eq:eta_even_moment_bound}
    \E \eta^{2k}
    \leq 2k \, \sigma^{2k} \int\limits_0^{+\infty} u^{k-1} e^{-u} \dd u
    = 2k \, \sigma^{2k} \, \Gamma(k)
    = 2 \cdot k! \cdot \sigma^{2k}.
\end{equation}

\medskip

\noindent
\textbf{Step 2: a bound on $\E \eta^{2}$.}
\quad
For the second moment, it is possible to prove a bit better bound, than \eqref{eq:eta_even_moment_bound}. It holds that
\[
    \E \eta^2
    = \int\limits_0^{+\infty} \p\left( \eta^{2} \geq t \right) \dd t
    \leq \int\limits_0^{\sigma^2 \log 2} \dd t
    + \int\limits_0^{\sigma^2 \log 2} \p\left( \eta^{2} \geq t \right) \dd t
    \leq \sigma^2 \log 2 + 2 \int\limits_{\sigma^2 \log 2}^{+\infty} e^{-t / \sigma^2} \dd t.
\]
Let us replace $t$ with $\sigma^2 (u + \log 2)$, where $u \in (0, + \infty)$. Then
\begin{equation}
    \label{eq:eta_2nd_moment_bound}
    \E \eta^2
    \leq \sigma^2 \log 2 + \sigma^2 \int\limits_0^{+\infty} e^{-u} \dd u
    = \sigma^2 (1 + \log 2).
\end{equation}

\medskip

\noindent
\textbf{Step 3: a bound on exponential moments.}
\quad
Let us fix any $\lambda > 0$ and consider $\E e^{\lambda \eta}$. Using Taylor's expansion and taking into account that $\E \eta = 0$, we obtain that 
\[
    \E e^{\lambda \eta}
    = 1 + \sum\limits_{k = 2}^\infty \frac{\lambda^k \E \eta^k}{k!}
    = 1 + \sum\limits_{k = 1}^\infty \frac{\lambda^{2k} \E \eta^{2k}}{(2k)!} + \sum\limits_{k = 1}^\infty \frac{\lambda^{2k + 1} \E \eta^{2k + 1}}{(2k + 1)!}.
\]
According to the Cauchy-Schwarz inequality, we have
\[
    \lambda^{2k + 1} \E \eta^{2k + 1}
    \leq \frac{\lambda^{2k} \E \eta^{2k}}2 + \frac{\lambda^{2k + 2} \E \eta^{2k + 2}}2,
\]
and thus,
\begin{align*}
    \E e^{\lambda \eta}
    &
    = 1 + \sum\limits_{k = 1}^\infty \frac{\lambda^{2k} \E \eta^{2k}}{(2k)!} + \sum\limits_{k = 1}^\infty \frac{\lambda^{2k + 1} \E \eta^{2k + 1}}{(2k + 1)!}
    \\&
    \leq 1 + \frac{7 \lambda^2 \E \eta^2}{12} + \sum\limits_{k = 2}^\infty \left( \frac1{(2k)!} + \frac{1}{2 (2k - 1)!} + \frac{1}{2 (2k + 1)!} \right)  \lambda^{2k} \, \E \eta^{2k}.
\end{align*}
Due to \eqref{eq:eta_even_moment_bound}, it holds that
\begin{align*}
    \E e^{\lambda \eta}
    &
    \leq 1 + \frac{7 \lambda^2 \E \eta^2}{12} + \sum\limits_{k = 2}^\infty \left( \frac1{(2k)!} + \frac{1}{2 (2k - 1)!} + \frac{1}{2 (2k + 1)!} \right)  \lambda^{2k} \cdot \left( 2 \cdot k! \cdot \sigma^{2k} \right)
    \\&
    \leq 1 + \frac{7 \lambda^2 \E \eta^2}{12} + \sum\limits_{k = 2}^\infty \frac{2 (k!)^2}{(2k)!} \left( 1 + k + \frac{1}{2 (2k + 1)} \right)  \cdot \frac{\lambda^{2k} \sigma^{2k}}{k!}.
\end{align*}
Note that, for any integer $k \geq 2$, we have
\begin{align*}
    \frac{2 (k!)^2}{(2k)!} \left( 1 + k + \frac{1}{2 (2k + 1)} \right)
    &
    = \left( 2 + 2k + \frac{1}{2k + 1} \right) \prod\limits_{j = 1}^k \frac{j \cdot j}{(2j - 1) \cdot 2j}
    \\&
    = \frac{1}{2^k} \left( 2 + 2k + \frac{1}{2k + 1} \right) \prod\limits_{j = 1}^k \frac{j }{(2j - 1)}
    \\&
    \leq \frac{2k + 3}{2^k} \prod\limits_{j = 1}^k \frac{j }{(2j - 1)}
    \leq \frac{2k + 3}{2^k} \cdot \frac13 < 1.
\end{align*}
Hence, it holds that
\begin{align*}
    \E e^{\lambda \eta}
    &
    \leq 1 + \frac{7 \lambda^2 \E \eta^2}{12} + \sum\limits_{k = 2}^\infty \left( \frac1{(2k)!} + \frac{1}{2 (2k - 1)!} + \frac{1}{2 (2k + 1)!} \right)  \lambda^{2k} \cdot \left( 2 \cdot k! \cdot \sigma^{2k} \right)
    \\&
    \leq 1 + \frac{7 \lambda^2 \E \eta^2}{12} + \sum\limits_{k = 2}^\infty \frac{\lambda^{2k} \sigma^{2k}}{k!}.
\end{align*}
The sharper bound \eqref{eq:eta_2nd_moment_bound} on $\E \eta^2$ implies that
\[
    \frac{7 \lambda^2 \E \eta^2}{12}
    \leq \frac{7 \lambda^2 \sigma^2 (1 + \log 2)}{12}
    \leq \lambda^2 \sigma^2.
\]
Thus, we finally obtain that
\[
    \E e^{\lambda \eta}
    \leq 1 + \lambda^2 \sigma^2 + \sum\limits_{k = 2}^\infty \frac{\lambda^{2k} \sigma^{2k}}{k!}
    \leq e^{\lambda^2 \sigma^2}
    \quad \text{for any $\lambda \geq 0$.}
\]
The proof for negative $\lambda$ is absolutely similar.

\endproof

\end{document}